\documentclass[final]{etds}

\usepackage{graphicx}   
\usepackage{euscript}  
\usepackage{amsfonts}  
\usepackage{amssymb}
\usepackage{latexsym}
\usepackage{amscd}
\newcommand{\field}[1]{\mathbb{#1}} 
\newtheorem{thm}{Theorem}[section] 
\newtheorem{nthm}{Theorem}

\newtheorem{lem}[nthm]{Lemma}

\newtheorem{defn}[nthm]{Definition}
 
\newtheorem{rmk}[nthm]{Remark}

\newtheorem{clm}{Main Claim}
\newtheorem{mainobs}{Numerical Observation}

\def\P{\tilde{P}}

\def\eps{\epsilon}

\def\Tble{T_{b,\lambda}}
\def\oT{\cT_{\tau}}
\def\Veut{V_{u,\tau}}

\def\gmu{\mu}
\def\mF{\mathfrak{F}}

\def\Um{\mathfrak{U}}
\def\mf{\mathfrak{f}}
\def\mc{\mathfrak{c}}
\def\um{\mathfrak{u}}

\def\mL{\mathfrak{L}}

\def\mB{\mathfrak{B}}
\def\mw{\mathfrak{w}}
\def\mW{\mathfrak{W}}
\def\mm{\mathfrak{m}}
\def\mM{\mathfrak{M}}
\def\mK{\mathfrak{K}}
\def\mE{\mathfrak{E}}
\def\mv{\mathfrak{v}}
\def\mS{\mathfrak{S}}
\def\mV{\mathfrak{V}}
\def\my{\mathfrak{y}}
\def\mY{\mathfrak{Y}}
\def\mDu{\mathfrak{Du}}
\def\mDf{\mathfrak{Df}}

\def\fC{\field{C}}

\def\fR{\field{R}}

\def\Ep{\mathcal{EP}}
\def\cA{\mathcal{A}}
\def\cB{\mathcal{B}}
\def\cC{\mathcal{C}}
\def\cD{\mathcal{D}}
\def\cE{\mathcal{E}}
\def\cF{\mathcal{F}}
\def\cG{\mathcal{G}}
\def\cH{\mathcal{H}}
\def\cK{\mathcal{K}}

\def\cL{\mathcal{L}}
\def\cN{\mathcal{N}}
\def\cM{\mathcal{M}}
\def\cT{\mathcal{T}}
\def\cR{\mathcal{R}}
\def\cRG{\mathcal{RG}}
\def\cP{\mathcal{P}}
\def\cS{\mathcal{S}}

\def\cW{\mathcal{W}}
\def\cV{\mathcal{V}}
\def\cX{\mathcal{X}}
\def\cZ{\mathcal{Z}}
\def\cO{\mathcal{O}}
\def\cQ{\mathcal{Q}}

\begin{document}

\ETDS{1}{1}{vol. num}{2008}

\large
\runningheads{D. Gaidashev, H. Koch}{Period Doubling in Area-Preserving Maps}
\title{Period Doubling in Area-Preserving Maps: An Associated One Dimensional Problem}
\author{Denis Gaidashev\affil{1}, Hans Koch\affil{2}}
\address{
\affilnum{1}Department of Mathematics,
University of Uppsala, Uppsala, Sweden.\\
 \email{gaidash@math.kth.se}\\
\affilnum{2}Department of Mathematics,
University of Texas at Austin, Austin, Texas, USA.\\
 \email{koch@math.utexas.edu}\\
}

\recd{}



\begin{abstract}

It has been observed that the famous Feigenbaum-Coullet-Tresser period doubling universality has a counterpart for area-preserving maps of ${\fR}^2$. A renormalization approach has been used in a computer-assisted proof of existence of an area-preserving map with orbits of all binary periods in \cite{EKW1} and \cite{EKW2}. As it is the case with all non-trivial universality problems in non-dissipative systems in dimensions more than one, no analytic proof of this period doubling universality exists to date. 

We argue that the period doubling renormalization fixed point for area-preserving maps is almost one dimensional, in the sense that it is close to the following H\'enon-like map:
$$H^*(x,u)=(\phi(x)-u,x-\phi(\phi(x)-u )),$$
where $\phi$ solves
$$\phi(x)={2 \over \lambda} \phi(\phi(\lambda x))-x.$$ 

We then give a ``proof'' of existence of solutions of small analytic perturbations of this one dimensional problem,  and describe some of the properties of  this solution. 
 The ``proof'' consists of an analytic argument for factorized inverse branches of $\phi$  together with verification of several inequalities and inclusions of subsets of $\field{C}$  numerically. 

Finally, we suggest an analytic approach to the full period doubling problem for area-preserving maps based on its proximity to the one dimensional.  In this respect, the paper is an exploration of a possible analytic machinery for a non-trivial renormalization problem in  a conservative two-dimensional system.

\end{abstract}

\newpage

\tableofcontents

\setcounter{page}{1}

\section{Introduction} 
Following the pioneering discovery of the Feigenbaum-Coullet-Tresser period doubling universality in unimodal maps \cite{Fei1}, \cite{Fei2}, \cite{TC} universality has been demonstrated to be a rather generic phenomenon in dynamics.

To prove universality one usually introduces a {\it renormalization} operator on a functional space, and demonstrates that this operator has a hyperbolic fixed point.

Such renormalization approach to universality has been very successful in one-dimensional dynamics, and has led to explanation of universality in unimodal maps \cite{Eps1}, \cite{Eps2},\cite{Lyu}, critical circle maps \cite{dF1,dF2}, \cite{Ya1}, \cite{Ya2} and holomorphic maps with a Siegel disks \cite{McM}, \cite{Ya3}, \cite{GaiYa}.

Universality has been abundantly observed in higher dimensions, in particular, in two and more dimensional dissipative systems (cf.  \cite{CEK1}, \cite{Spa}), in area-preserving maps, both as the period-doubling universality \cite{DP}, \cite{Hel}, \cite{BCGG}, \cite{CEK2}, \cite{EKW1}, \cite{EKW2}, and as the universality associated with the break-up of invariant surfaces  \cite{Shen}, \cite{McK1}, \cite{McK2}, \cite{ME},  and in Hamiltonian flows \cite{ED},\cite{AK}, \cite{AKW}, \cite{Koch1}, \cite{Koch2}, \cite{Koch3}, \cite{GK}, \cite{Gai1}, \cite{Kocic} . It has been established that the universal behavior in dissipative and conservative higher dimensional systems is fundamentally different. The case of of the dissipative systems is often reducible to the one-dimensional Feigenbaum-Coullet-Tresser universality (\cite{CEK1}, \cite{dCLM} ). The latter case is very different, and at present there is no deep understanding of universality in conservative systems, other than in the ``trivial'' case of the universality for systems ``near integrability'' \cite{Koch1}, \cite{Koch2}, \cite{Gai1}, \cite{Kocic}, \cite{KLDM} . The study of the interesting cases of universality for maps and flows far from linear, is at present confined to numerics (for instance, \cite{McK1},\cite{CEK1}, \cite{GK})), or computer-assisted proofs \cite{EKW2}, \cite{Koch2}, \cite{Koch3}. The latter approach turned to be quite powerful, albeit rather specialized and time and effort consuming. Very little analytic machinery exists for higher-dimensional maps and flows far from linear.

In this paper we will consider a period-doubling universality for area-preserving maps of the plane ---  an analogue of Feigenbaum-Coullet-Tresser universality in higher dimensions. 

An infinite period-doubling cascade in families of area-preserving maps was observed by several authors in  early 80's \cite{DP}, \cite{Hel}, \cite{BCGG}, \cite{Bou}, \cite{CEK2}. A typical period-doubling scenario can be illustrated with the area-preserving H\' enon family (cf. \cite{Bou}) :
$$ H_a(x,y)=(-y +1 - a x^2, x).$$

Maps in this family posses a fixed point $((-1+\sqrt{1+a})/a,(-1+\sqrt{1+a})/a) $ which is stable for $-1 < a < 3$. When $a_1=3$ this fixed point becomes unstable, at the same time an orbit of period two is born with $H_a(x_\pm,x_\mp)=(x_\mp,x_\pm)$, $x_\pm= (1\pm \sqrt{a-3})/a$. This orbit, in turn, becomes unstable at $a_2=4$, giving birth to a period $4$ stable orbit. Generally, there  exists a sequence of parameter values $a_k$, at which the orbit of period $2^{k-1}$ turns unstable, while at the same time a stable orbit of period $2^k$ is born. The parameter values $a_k$ accumulate on some $a_\infty$. The crucial observation is that the accumulation rate
\begin{equation}
\lim_{k \rightarrow \infty}{a_k-a_{k-1} \over  a_{k+1}-a_k } = 8.721...
\end{equation} 
is universal for a large class of families, not necessarily H\'enon.

Furthermore, the $2^k$ periodic orbits scale asymptotically with two scaling parameters
\begin{equation}
\lambda=-0.249...,\quad \gmu=0.061...
\end{equation}

\begin{figure}
\vspace{-1.0cm}
 \begin{center}
\begin{tabular}{c c c}
$\!\!\!\!\!\!\!\! $ \resizebox{45mm}{!}{\includegraphics{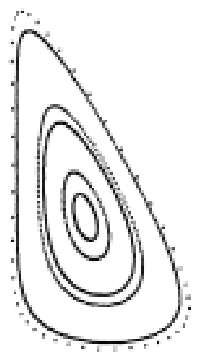}}$\!\!\!\!\!\!\!\! $&$\!\!\!\!\!\!\!\!$ \resizebox{45mm}{!}{\includegraphics{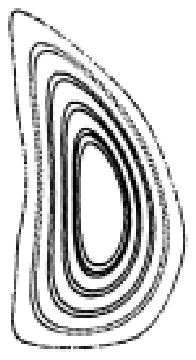}}$ \!\!\!\!\!\!\!\! $&$\!\!\!\!\!\!\!\! $\resizebox{45mm}{!}{\includegraphics{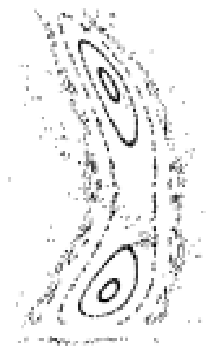}} 
\end{tabular} 
\vspace{-1.0cm}
\caption{\it Bifurcation of a stable fixed point into a stable period $2$ orbit  in the area-preserving H\'enon family $H_a$.}
\end{center}
\end{figure}

To explain how orbits scale with $\lambda$ and $\mu$ we will follow \cite{Bou}. Consider an interval $(a_k,a_{k+1})$ of parameter values in a `''typical'' family $F_a$. For any value $\alpha \in (a_k,a_{k+1})$ the map $F_\alpha$ posses a stable periodic orbit of period $2^k$. We fix some $\alpha_k$ within the interval $(a_k,a_{k+1})$ in some consistent way; for instance, by requiring that the restriction of $F^{2^k}_{\alpha_k}$ to a neighborhood of a stable periodic point in the $2^{k}$-periodic orbit is conjugate via a diffeomorphism $H_k$ to a rotation with some fixed rotation number $r$.  Let $p'_k$ be some unstable periodic point in the $2^{k-1}$-periodic orbit, and let $p_k$ be the further of the two stable $2^{k}$-periodic points that bifurcated from $p'_k$.  Then,
$${1 \over \lambda}=-\lim_{k \rightarrow \infty}{ |p_k-p'_k | \over |p_{k+1}-p'_{k+1} |},\quad {\lambda \over \mu}=-\lim_{k \rightarrow \infty}{ \rho_k \over \rho_{k+1}},$$
where $\rho_k$ is the ratio of the eigenvalues of $D H_k(p_k)$.

This universality can be explained rigorously if one shows that the {\it renormalization} operator
\begin{equation}\label{Ren}
R[F]=\Lambda^{-1}_F \circ F \circ F \circ \Lambda_F,
\end{equation}
where $\Lambda_F$ is some $F$-dependent coordinate transformation, has a fixed point, and the derivative of this operator is hyperbolic at this fixed point.

It has been argued in \cite{CEK2}  that $\Lambda_F$ is a diagonal linear transformation. Furthermore, such $\Lambda_F$ has been used in \cite{EKW1} and \cite{EKW2} in a computer assisted proof of existence of the renormalization fixed point. This is the strongest result to date.

In this paper we will first present some numerical and qualitative evidence that the fixed point for the renormalization operator $(\ref{Ren})$ is, in the some appropriate sense, very close to an area-preserving H\'enon-like map
\begin{equation}
H^*(x,u)=(\phi(x)-u,x-\phi(\phi(x)-u )),
\end{equation}
where $\phi$ solves the following one-dimensional problem of non-Feigenbaum type:
\begin{equation}\label{0_equation}
\phi(y)={2 \over \lambda} \phi(\phi(\lambda y)) -y.
\end{equation}

Furthermore, we will consider a more general functional equation
\begin{equation}\label{tau_equation}
\phi(y)={2 \over \lambda} \phi(\phi(\lambda y)) - y+\tau(y),
\end{equation}
where $\tau$ is some small analytic ``perturbation''.

The problem $(\ref{tau_equation})$ will be reformulated as a fixed point problem for the diffeomorphic part of the factorized inverse branches of $\phi$. Specifically, we will restate the problem as a fixed point problem for a continuous operator on the diffeomorphic part, and demonstrate that there is a choice of a compact functional space which is left invariant by this operator. Existence of a fixed point will follow from the Schauder-Tikhonov Theorem.  Thus, the proof  relies on an analytic argument, however several rather technical conditions will be verified on a computer numerically.

Furthermore, we will suggest an approach to the full problem $R[F]=F$ based on a ``factorization'' of $F$ into a dominant ``one-dimensional'' H\'enon-like part plus a small two-dimensional correction. Qualitatively, the corrective term is contracted by the renormalization operator on a neighborhood of the ``one-dimensional'' approximate fixed point. An important remark is that this contraction is not shown by brute force (as that would be the case with a ``hard'' computer-assisted proof); instead, the period-doubling operator essentially turns into a contraction after the above-mentioned one-dimension problem is ``factored out''.

We will briefly outline further steps that are necessary to complete the proof of existence of a period-doubling fixed point. The completion of the ``program'' will is a subject of our upcoming work.

We should note that construction of an approximate renormalization fixed point is generally an easy task. What we are trying to do here is i) construct an approximate fixed point which lies in a compact functional space, ii) turn the period doubling operator for area-preserving maps into a contraction using the {\it a-priori} bounds from this compact space. The dominant ``one-dimensional'' part of the fixed point of the full problem $(\ref{Ren})$ is expected to lie in the same space and obey similar {\it a-priori} bounds.


\section{A renormalization operator for generating functions} 

An ``area-preserving map'' will mean an exact symplectic diffeomorphism of a subset of ${\fR}^2$ onto its image.

Recall, that an area-preserving map can be uniquely specified by its generating function $S$:
\begin{equation}\label{gen_func}
\left( x \atop -S_1(x,y) \right) {{ \mbox{{\small \it  F}} \atop \mapsto} \atop \phantom{\mbox{\tiny .}}} \left( y \atop S_2(x,y) \right), \quad S_i \equiv \partial_i S.
\end{equation}

Furthermore, we will assume that $F$ is reversible, that is 
$$T \circ F \circ T=F^{-1}, \quad {\rm where} \quad T(x,u)=(x,-u).$$

For such maps it follows from $(\ref{gen_func})$ that 
$$S_1(y,x)=S_2(x,y) \equiv s(x,y),$$
and
$$
\nonumber \left({x  \atop  -s(y,x)} \right)  {{ \mbox{{\small \it  F}} \atop \mapsto} \atop \phantom{\mbox{\tiny .}}} \left({y \atop s(x,y) }\right).
$$

It is this ``little'' $s$ that will be referred to below as ``the generating function''.

We will now derive an equation for the generating function of the renormalized map $\Lambda_F \circ F \circ F \circ \Lambda_F^{-1}$.

Applying a reversible $F$ twice we get
$$
 \left({x'  \atop  -s(z',x')} \right) {{ \mbox{{\small \it  F}} \atop \mapsto} \atop \phantom{\mbox{\tiny .}}} \left({z' \atop s(x',z')} \right)=
\left({z'  \atop  -s(y',z')} \right) {{ \mbox{{\small \it  F}} \atop \mapsto} \atop \phantom{\mbox{\tiny .}}} \left({y'\atop  s(z',y')} \right).
$$

It has been argued in \cite{CEK2}  that 

$$\Lambda(x,u)=(\lambda x, \gmu u).$$

We therefore set  $(x',y')=(\lambda x,  \lambda y)$, $z'(\lambda x, \lambda y)= z(x,y)$ to obtain:

\begin{equation}\label{doubling}
\left(\!{x  \atop  -{ 1 \over \gmu } s(z,\lambda x)} \!\right) \!{{ \mbox{{\small $\Lambda$}} \atop \mapsto} \atop \phantom{\mbox{\tiny .}}} \!\left(\!{\lambda x  \atop  -s(z,\lambda x)} \!\right) \!{{ \mbox{{\small \it  F $ \circ$ F}} \atop \mapsto} \atop \phantom{\mbox{\tiny .}}}\!\left(\!{\lambda y \atop s(z,\lambda y)}\! \right)   {{ \mbox{{\small \it  $\Lambda^{-1}$}} \atop \mapsto} \atop \phantom{\mbox{\tiny .}}} \left(\!{y \atop {1 \over \gmu } s(z,\lambda y) }\!\right),
\end{equation}
where $z(x,y)$ solves
\begin{equation}\label{midpoint}
s(\lambda x, z(x,y))+s(\lambda y, z(x,y))=0.
\end{equation}

If the solution of $(\ref{midpoint})$ is unique, then $z(x,y)=z(y,x)$, and it follows from $(\ref{doubling})$ that the generating function of the renormalized $F$ is given by 
\begin{equation}
\tilde{s}(x,y)=\gmu^{-1} s(z(x,y),\lambda y).
\end{equation}

Furthermore, it is possible to fix some normalization conditions for $\tilde{s}$ and $z$ which serve to determine scalings $\lambda$ and $\gmu$ as functions of $s$. Notice, that the normalization
$$s(1,0)=0$$ 
is reproduced for $\tilde{s}$ as long as 
$$z(1,0)=z(0,1)=1.$$

In particular, this implies that 
$$s(\lambda, 1)+s(0, 1)=0.$$

Furthermore, the condition
\begin{equation}\label{cond_s1}
\partial_1 s(1,0)=1
 \end{equation}
is reproduced as long as 
$$\gmu=\partial_1 z (1,0).$$

We will now summarize the above discussion in the following definition of the renormalization operator acting on generating functions originally due  to the authors of \cite{EKW1} and \cite{EKW2}:
\begin{defn}\begin{eqnarray}
\label{ren_eq} {\cR}_{EKW}[s](x,y)&=&\gmu^{-1} s(z(x,y),\lambda y), {\rm where}\cr
\label{midpoint_eq} 0&=&s(\lambda x, z(x,y))+s(\lambda y, z(x,y)), \cr
0&=&s(\lambda,1)+s(0,1) \quad {\rm and} \quad \gmu=\partial_1 z (1,0).
\end{eqnarray}
\end{defn}

As we have already mentioned the following has been proved with the help of a computer in \cite{EKW1} and \cite{EKW2}:
\begin{thm}
There is an $s^*$ in some Banach space of analytic functions, such that the operator ${\cR}_{EKW}$ is well-defined, analytic and compact on some neighborhood of $s^*$, and ${\cR}_{EKW}[s^*]=s^*.$ Furthermore, the scalings $\lambda^*$ and $\gmu^*$ corresponding to the fixed point $s^*$ satisfy
\begin{eqnarray}
\label{lambda} -0.2492 < &\lambda& < -0.2485, \\
\label{mu} 0.0606 < &\gmu& < 0.0616.
\end{eqnarray}
  \end{thm}

Here, we will quote for reference purposes approximations of several first coefficients of the fixed point $s^*$ and the corresponding midpoint function $z^*$, together with the set off relations between them, all of which are obtained by differentiation of the fixed point equation $s(x,y)=\gmu^{-1} s(z(x,y),\lambda y)$:
\begin{eqnarray}
 \label{coeffs1} s^*(x,y)\!\!&\!\!=\!\!&\!\!(x\!-\!1)+a y+ {b \over 2} (x\!-\!1)^2+c(x\!-\!1)y +{d \over 2} y^2 +O((x\!-\!1)^i y^j), \!\quad \! j\!+\!i=\!3,\\
 \label{coeffs2} z^*(x,y)\!\!&\!\!=\!\!&\!\!1+\gmu(x\!-\!1)+\theta y + {\upsilon \over 2}(x\!-\!1)^2 + \iota(x\!-\!1) y + {\nu \over 2} y^2+O((x\!-\!1)^i y^j),\! \! \quad \!\!j\!+\!i\!=\!3,
\end{eqnarray}
\begin{displaymath}
\begin{array}{lllllllll}
 &&a=0.1948..., \quad   &&b=-0.0523..., \quad &&c=-0.0497, \quad &&  d=2.11...,\\
 &&a=  {\theta \over \gmu-\lambda},  \quad  &&b \gmu={\upsilon \over 1-\gmu},  \quad  &&c \gmu ={\upsilon \theta+\iota(1-\gmu) \over (1-\lambda) (1-\gmu)}, \quad  &&d=a(2 c -b a)+{2 \iota a-\upsilon a^2-\nu \over \lambda^2-\gmu }.
\end{array}
\end{displaymath}


\section{A $\lambda$-manifold based approach. Nonlinear scalings.}\label{manifold}

The fixed point equation ${\cR}_{EKW}[s]=s$ together with the tautological identity $\lambda y=\lambda y$ can be written simultaneously as
\begin{equation}\label{Schroder_eq}
A S=S \circ G, \quad A=\left[{\gmu \quad  0 \atop 0  \quad \lambda} \right],
\end{equation}
where
$$G(x,y)=(z(x,y),\lambda y) \quad {\rm and} \quad S(x,y)=(s(x,y),y).$$

Given $G$ this is a  equation for a linearizer $S$ that conjugates $G$ to a linear map $A$. Notice that the point $(1,0)$ is fixed under $G$, the line $y=0$ is, locally, an invariant manifold for $G$ associated with the eigenvalue $\gmu=\partial_1 z (1,0)$. Next, suppose that the function $\phi(y)$ defined on some neighborhood of $0$, is such that $s(\phi(y),y)=0$ (necessarily, $\phi(0)=1$). Then, the first equation of $(\ref{Schroder_eq})$, evaluated at $(\phi(y),y)$ becomes
$$s(G(\phi(y),y))=s(z(\phi(y),y),\lambda y)=0,$$
which implies that
$$z(\phi(y),y)=\phi(\lambda y).$$
That is 
$$G(\phi(y),y)=(\phi(\lambda y),\lambda y),$$
and the curve $(\phi(y),y)$ is an invariant manifold for $G$ associated with the eigenvalue $\lambda$. 

We will refer to the two invariant manifolds as the $\lambda$- and the $\gmu$-manifolds.

A solution $S$ of the Schr\"oder equation $(\ref{Schroder_eq})$ maps the fixed point of $G$ to the origin, the $\gmu$-manifold to the $x$-axis, the $\lambda$-manifold --- to the $y$-axis. Such (normalized) solution necessarily satisfies 
\begin{equation}\label{s}
s(x,y)=(x-\phi(y))(1+\varepsilon(x,y)) \equiv x-\phi(y) + \epsilon(x,y), 
\end{equation}
where $\phi$ is the parametrization of the $\lambda$-manifold of $G$.

Notice that the fixed point problem for the renormalization operator $(\ref{ren_eq})$ is equivalent to the solution of the Schr\"oder equation $(\ref{Schroder_eq})$ together with the midpoint equation $(\ref{midpoint_eq})$.

The second of these equations can be readily solved on some neighborhood of $(1,0)$ for a rather specific $s(x,y)=x-\phi(y)$ ($\epsilon(x,y)=0$), where $\phi$ is some function, invertible on a neighborhood of zero (with no further assumptions on $\phi$ at this point). The midpoint equation for such $s$ becomes
$$
x+y={2 \over \lambda} \phi(z(x,y)).
$$ 

Evaluate this midpoint equation at $(\phi(y),y)$:
$$
\phi(y)+y={2 \over \lambda} \phi(z(\phi(y),y)).
$$

{\it Clearly,  if $\phi(y)$ is also the parametrization of the $\lambda$-manifold  for $G(x,y)=(z(x,y),\lambda y)$, then it satisfies the following functional equation:}
\begin{equation}\label{func_eq}
\phi(y)={2 \over \lambda} \phi(\phi(\lambda y))-y.
\end{equation}

{\it We would like to reiterate that a solution $\phi$ of $(\ref{func_eq})$} {\it is the parametrization of the $\lambda$-manifold for $G(x,y)=(z(x,y),\lambda)$ where $z$ solve the midpoint equation for $s(x,y)=x-\phi(y)$, but such $s$ is not yet the solution of the Schr\"oder equation, nor of the fixed point problem for ${\cR}$}.

Now, suppose that $\epsilon(x,y)$ in  $(\ref{s})$ is non-zero. Then the midpoint equation for such $s$  reads:
$$
x+y={2 \over \lambda} \phi(z)-{1 \over \lambda}\epsilon(\lambda x,z) -{1 \over \lambda} \epsilon(\lambda y,z).
$$ 

Again, suppose that $\phi$ is also a parametrization for the $\lambda$-manifold of $G$, then it satisfies:

\begin{equation}\label{new_eq}
\phi(y)={2 \over \lambda } \phi(\phi(\lambda y))-y-\omega_{\eps,\phi}(y),
\end{equation}
where
\begin{eqnarray}
\nonumber  \omega_{\eps,\phi}(y) & \equiv& {1\over \lambda} \epsilon(\lambda \phi(y),\phi(\lambda y) )+ {1 \over \lambda} \epsilon(\lambda y, \phi(\lambda y) ).
\end{eqnarray}


At this point we will modify the operator ${\cR}_{EKW}$ by introducing a nonlinear scaling $\gmu=\gmu(y)$ and turning it into an operator for the ``corrective term'' $\epsilon$. The reasons for this modification will be clear momentarily.
\begin{defn}\label{ren_op}
Assume that given $\epsilon$ there is a solution $\varphi_\eps$ of the equation $(\ref{new_eq})$ defined on some neighborhood of the interval $(0,1)$. Set, formally,
\begin{equation}\label{RG}
{\cRG}[\epsilon](x,y) \equiv  \mu^{-1}(y)(z(x,y)-\varphi_\eps(\lambda y)+\epsilon(z(x,y),\lambda y)) -(x-\varphi_\eps(y)),
\end{equation}
where the midpoint function $z$ solves
\begin{equation}\label{new_midpoint}
\lambda x -\varphi_\eps(z(x,y))+ \lambda y -\varphi_\eps(z(x,y)) +\eps(\lambda x, z(x,y))+ \eps(\lambda y,z(x,y))=0,
\end{equation}
and  scalings $\lambda$ and $\gmu(y)$  satisfy 
\begin{eqnarray}
\lambda&=&2 \varphi_\eps(1)-\eps(\lambda,1)-\eps(0,1),\\
\gmu(y)&=&\partial_1 z (\phi(y),y).
\end{eqnarray}
\end{defn}

Suppose that there is a fixed point $\epsilon^*$ such that  ${\cRG}[\epsilon^*]=\epsilon^*$, the corresponding scaling being $\lambda^*$ and $\gmu^*$, the midpoint function --- $z^*$. Then it is immediate that the function
$$s^*(x,y)=x-\phi_{\epsilon^*}(y)+\epsilon^*(x,y)$$
satisfies
\begin{equation}
s^*={s^* \circ G^* \over \gmu^*}, \quad {\rm where}\quad G^*(x,y) \equiv (z^*(x,y),\lambda^* y),
\end{equation}
and the corresponding reversible $F^*$ satisfies
\begin{eqnarray}
 F^*&=& \Lambda^{-1}_{*} \circ F^* \circ F^* \circ \Lambda_{*} \\
\Lambda_*(x,u)&=&(\lambda^* x, \mu^*(x) u), \\
\Lambda^{-1}_*(y,v)&=&\left( {y \over \lambda^*}, {v \over \mu^* \left( {y \over \lambda^*} \right) } \right).
\end{eqnarray}

A this point the Definition $\ref{ren_op}$ is purely formal. The standard properties of being well-defined and analytic (on a neighborhood of the fixed point) should be verified after an functional space for $\epsilon$'s is appropriately chosen.  

Notice, that $s(x,y)=x-\varphi_\eps(y) + \eps(x,y)$ is an approximate renormalization fixed point on a neighborhood of  the $\lambda$-manifold $(\varphi_\eps(y),y)$ for a any sufficiently small $\epsilon$, in the sense that
\begin{eqnarray}
 \nonumber \cR_{EKW}[s](x,y)&=&\mu^{-1}(y) (z(x,y)-\varphi_\eps(\lambda y)+\epsilon(z(x,y),\lambda y))\\
\nonumber &=&\mu^{-1}(y) (z (x ,y)-z(\varphi_\eps(y),y) +\epsilon(z(x,y) ,\lambda y))\\
\nonumber &=&{ \partial_1 z (\varphi_\eps(y),y) \over \mu(y)} (x-\varphi_\eps(y)+O((x-\varphi_\eps(y))^2) +\mu^{-1}(y) \eps(z(x,y),\lambda y) \\
\nonumber &=&x-\varphi_\eps(y)+O((x-\varphi_\eps(y))^2) +\mu^{-1}(y) \eps(z(x,y),\lambda y)). 
\end{eqnarray}

We would like to emphasize that the coefficients of $O$ are proportional to third and higher-order derivatives of $z$ which are ``tiny'' (cf. $(\ref{coeffs2})$). Therefore, $O$ is expected to be small in any reasonable norm on a sufficiently small neighborhood of the $\lambda$-manifold $(\varphi_\eps(y),y)$.

Furthermore, 
the operator $\cRG$ is expected to be a contraction on a neighborhood of the $\lambda$-manifold. This can be seen if one writes 
\begin{eqnarray}
\nonumber \epsilon(x,y) & \equiv & (x-\varphi_\eps(y)) \varepsilon(x,y),\\
\nonumber \varepsilon(x,y) &  = & f(y)+O((x-\varphi_\eps(y))),\\
\nonumber s(x,y)&=&(x-\varphi_\eps(y)) (1+f(y))+O( (x-\varphi_\eps(x,y))^2 ),
\end{eqnarray}
then the condition $(\ref{cond_s1})$ implies
\begin{equation}
1+f(0)=1=> f(y)=O(y) => \varepsilon(x,y)=O((x-\varphi_\eps(y))^i y^j), \quad i+j=1,
\end{equation}
in particular, 
\begin{equation}
\varepsilon(1,0)=0,
\end{equation}
and the operator 
$$ \varepsilon(x,y) \mapsto \varepsilon(G(x,y))$$  
is a contraction. Then,
\begin{eqnarray}
\nonumber \cRG[\epsilon](x,y)&=& O( (x-\varphi_\eps(y))^2 ) +\mu^{-1}(y) \eps(z(x,y),\lambda y))\\
\nonumber&=&O((x-\varphi_\eps(y))^2)+\mu^{-1}(y)\left(z(x,y)-\varphi_\eps(\lambda y) \right) \varepsilon(G(x,y))\\
\nonumber&=& O((x-\varphi_\eps(y))^2)+ { \partial_1 z (\varphi_\eps(y),y) \over \mu(y)} (x-\varphi_\eps(y)+O((x-\varphi_\eps(y))^2)) \varepsilon(G(x,y))\\
\nonumber&=&(x-\varphi_\eps(y)) \varepsilon(G(x,y)) + O((x-\varphi_\eps(y))^2),
\end{eqnarray}
is also a contraction on a sufficiently small neighborhood of the curve $(\varphi_\eps(y),y)$.




Any rigorous formulation of the above informal discussion should start with, first, solving the equation $(\ref{new_eq})$ and, second, identifying a suitable functional space for $\epsilon$'s. In this paper we will concentrate on studying the equation $(\ref{new_eq})$ and its special case $(\ref{func_eq})$. In Section ~\ref{towards}, we will suggest an iterative procedure for the operator $\cRG$ which can be used to prove the existence of a fixed point for this operator. Its implementation will be a subject of an upcoming work.


\section{Notation. Some facts about Herglotz-Pick functions}

We will proceed with some definitions.

We will use the standard notation for the lower and upper-half planes:
$$\fC_\pm \equiv \{z \in \fC: \pm {\Im (z)} > 0 \}.$$

Let $J=(l,r) \subset \fR$. Define $D_+(J,\theta)$ to be an open subset of $\fC_+$  bounded by a circular arc intersecting $\fR$ at the endpoints of $J$ at an angle $\theta$, and let $D_-(J,\theta)=D_+(J,\theta)^*$ where ${}^*$ stands for the complex conjugation.  A Poincar\'e neighborhood is defined as 
$$D(J,\theta)=D_+(J,\theta) \cup D_-(J,\theta) \cup J.$$

Given an interval $J \subset \fR$, denote
$$
\fC(J) \equiv \fC_+ \cup \fC_- \cup J, \quad \fC_1 \equiv \fC((-1,1)).   
$$ 

We will denote ${\cF}(\cD)$ the Frechet space of functions holomorphic on a domain $\cD$ equipped with the topology of uniform convergence on compacts. A subset of functions in $\cF$ assuming their values in a set $\cE$, will be denoted by $\cO(\cD,\cE)$.

Suppose that $\cD$ is real symmetric, and let ${\bf \mc}=\{\mc_1,\mc_2,\mc_3,\mc_4\}$ be a quadruple of real numbers, such that $\{\mc_1,\mc_2\} \in \cD$ and $\{\mc_3,\mc_4\} \in \cE$. We will further define

\begin{eqnarray}
\nonumber {\cA}(\cD,\cE;{\bf \mc}) &\equiv& \left\{u \in {\cO}(\cD,\cE): u(z)=u(z^*)^*, u(\cD \cap \fC_\pm) \subset \overline{\cE \cap \fC_\pm}, u \left( \mc_1 \right)=\mc_3, u(\mc_2)=\mc_4 \right\},\\
\nonumber  \cA_1({\bf c}) &\equiv& \cA (\fC_1,\fC_1;{\bf c}).
\end{eqnarray}

 Clearly, $\cA(\cD,\cE,{\bf \mc})$ is isomorphic to some $\cA_1({\bf c} )$ through a unique conformal isomorphism $\Phi$ that is normalized so that 
$$\Phi(l)=-1, \quad \Phi(r)=1, \quad \Phi(a)=b$$  
Here $a$ and $b$ are some constants, and
$$c_i=\Phi(\mc_i), \quad i=1..4.$$ 

Functions in ${\cA}_1({\bf c})$, admit the following integral representation:
\begin{equation}\label{int_rep}
f(z)-c_3= a (z-c_1) + \int d \nu(t) \left({1 \over t-z} - {1 \over t-c_1}   \right), 
\end{equation}
where $\nu$ is a measure supported in ${\fR} \setminus (-1,1)$.

This integral representation can be used to obtain the following {\it a-priori} bounds for functions in ${\cA}_1({\bf c})$
\begin{eqnarray} 
\label{function_1}  {c_4-c_3  \over c_2-c_1} {1+c_2 \over 1+x} \le &{f(x)-c_3 \over x-c_1} &\le  {c_4-c_3  \over c_2-c_1} {1-c_2 \over 1-x}, \quad x \in (-1,c_2),\\
 \label{function_2}  {c_4-c_3  \over c_2-c_1} {1+c_2 \over 1+x} \ge &{f(x)-c_3 \over x-c_1} &\ge  {c_4-c_3  \over c_2-c_1} {1-c_2 \over 1-x}, \quad x \in (c_2,1),\\
\label{first_der} { 1+c_1 \over (x-c_1 )  (1+x) } \le &{f'(x) \over f(x)-c_3}& \le { 1-c_1 \over (x-c_1 )  (1-x) }  , \quad x \in (-1,1),\\
\label{second_der} {-2 f'(x) \over 1+x}  \le &f''(x)& \le {2 f'(x) \over 1-x},  \quad x \in (-1,1).
\end{eqnarray}

If $\Phi \arrowvert_\fR$ is a monotone function, then one can transfer the bounds $(\ref{function_1})$---$(\ref{second_der})$ to $\cA(\cD,\cE,{\bf \mc})$.

Finally, we will mention the following version of Schwarz Lemma  which will play an important role in our results below (cf \cite{Eps2}, \cite{Sul}, \cite{LY}):

\begin{lem}\label{Epstein_lemma}
Let $u:\fC_J \mapsto \fC_{J'}$ be a holomorphic  map such that $u(J) \subset J'$. Then for any $\theta \in (0,\pi)$,  $u(D_\pm(J,\theta)) \subset D_\pm(J',\theta)$.
\end{lem}


\section{Main results and observations}

We will now summarize the main results of the paper.

Define
\begin{equation}
\label{K_op}{\cK}_\epsilon[\phi](y)={2 \over \lambda } \phi(\phi(\lambda y))-y-\omega_{\eps,\phi}(y),
\end{equation}
where $\omega_\phi$ is as in $(\ref{new_eq})$.

Below we will denote the solution of ${\cK}_\eps[\phi]=\phi$ by $\varphi_\eps$,  and that of equation $(\ref{tau_equation})$ by $\phi_\tau$.



\begin{figure}
 \begin{center}
\begin{tabular}{c c c}
 \resizebox{40mm}{!}{\includegraphics[angle=-90]{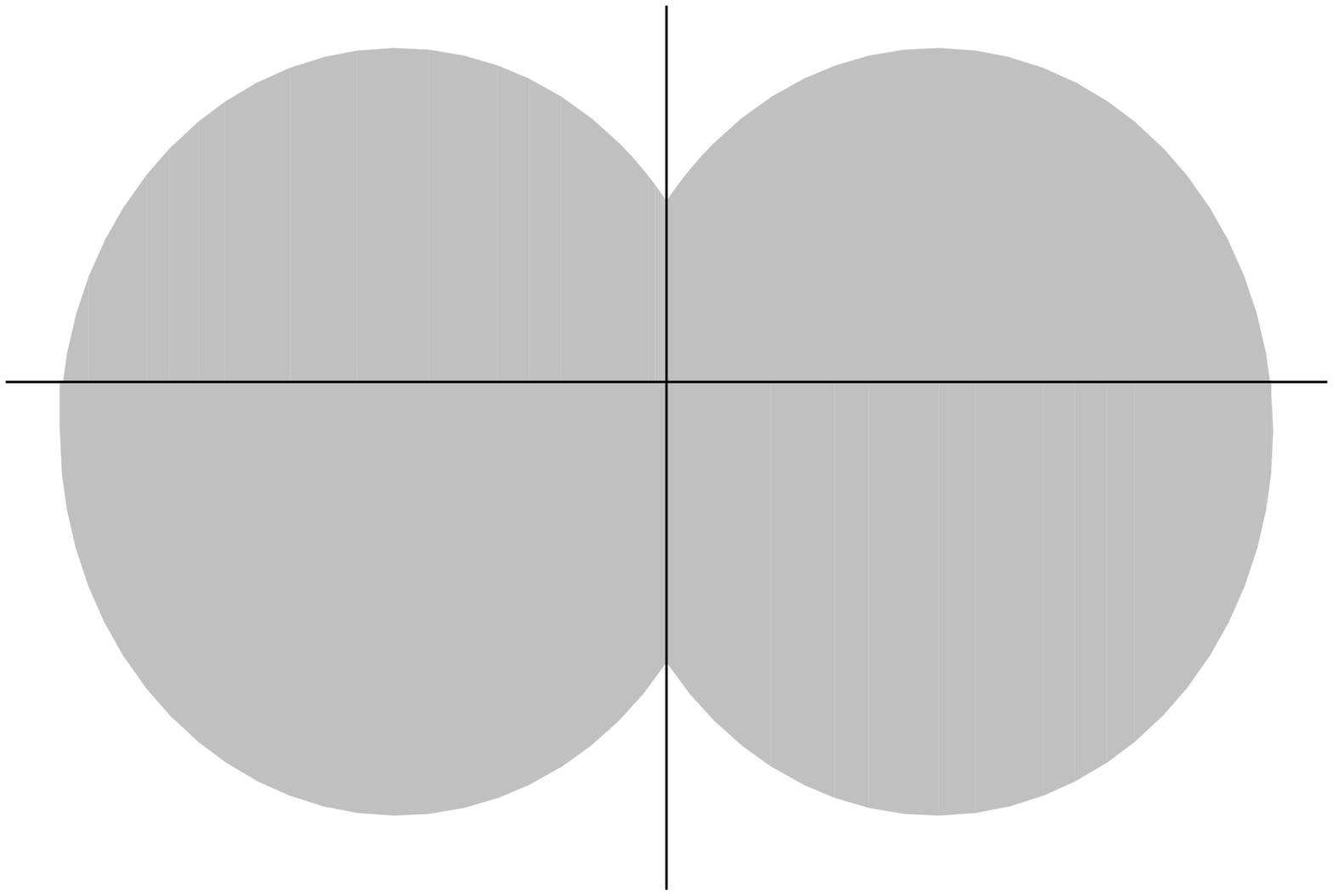}} & \phantom{aaaaaa}  &  \resizebox{35mm}{!}{\includegraphics[angle=-90]{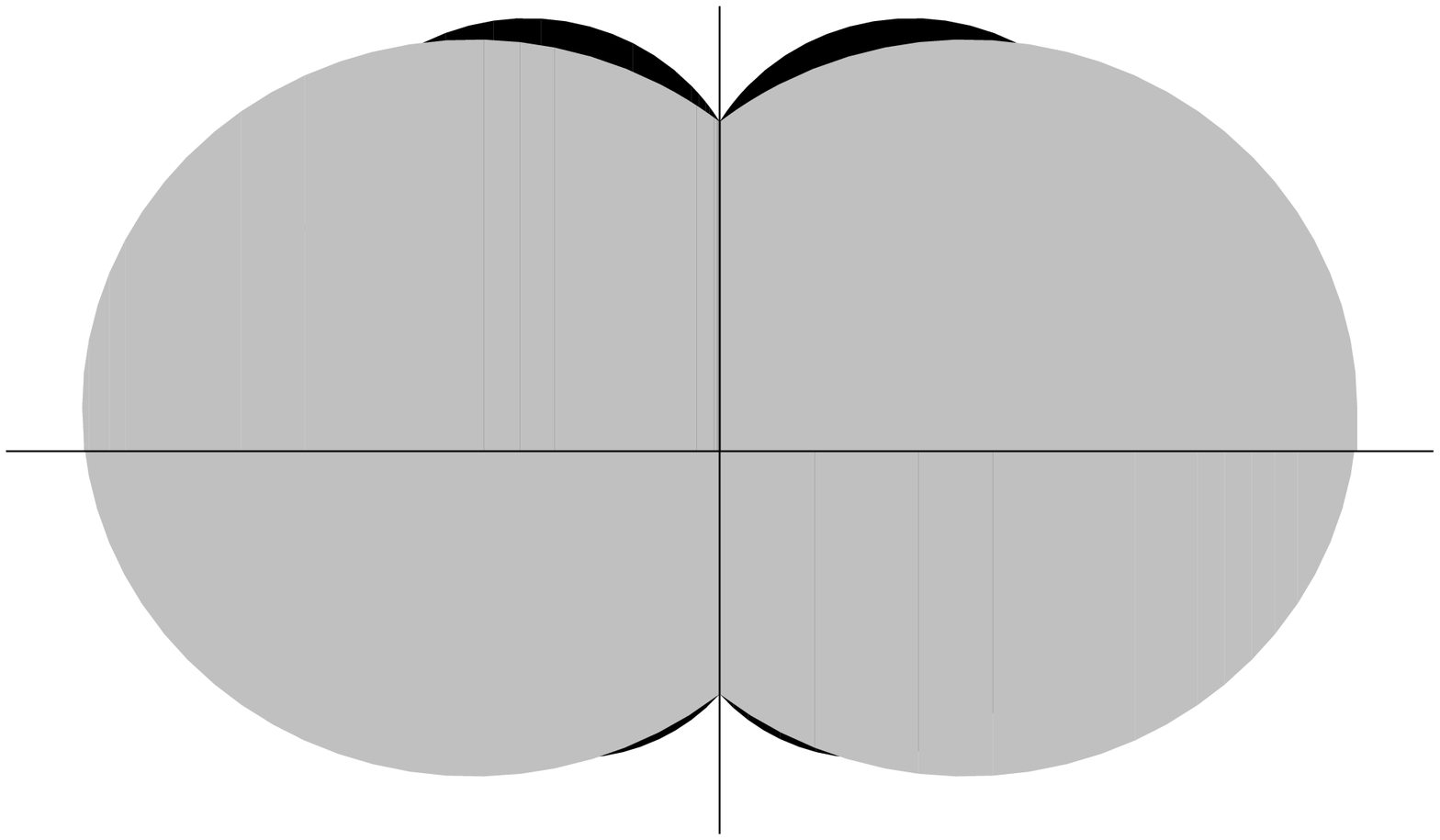}} \\
a) & &  b)
\end{tabular} 
\caption{\it Domain a) and range b) of functions in $\cA(\cD,\cE,{\bf \mc})$. In b) the Poincar\'e neighborhood $D(I_2,\theta_2)$ is given in light color, $D(I_3,\theta_3) \setminus D(I_2,\theta_2)$ and  $D(I_4,\theta_4) \setminus D(I_2,\theta_2)$ --- in dark.}
\end{center}
\end{figure}

\begin{clm}\label{main_thm}

Set
\begin{eqnarray}
\label{I1}I_1&=&(-1.49,0.96), \quad \theta_1= 0.2075 \pi,\\
\label{I2} I_2&=&(-2.347360978,3.181216988), \quad \theta_2= 0.28285  \pi,\\
\label{I3}I_3&=& (-2.347360978,0.652639022) , \quad \theta_3= 0.8,\\
\label{I4}I_4&=& (1.681216988,3.181216988) , \quad \theta_3= 0.4,
\end{eqnarray}
and $\cD=D1(I_1,\theta_1)$,  $\cE=D(I_2,\theta_2)  \cup D(I_3,\theta_3) \cup  D(I_4,\theta_4)$.

There exists $ \delta >0$ and $\kappa>0$, such that given  any $\tau$ holomorphic in $\cE$, real-valued on $\cE \cap \fR$ and satisfying  
$$\sup_{z \in \cE} |\tau(z)| < \delta, \quad \sup_{z \in \cE} |\tau'(z)| \le \kappa, \quad \tau(0)=0,$$
there exists a function $\phi_\tau$, holomorphic on some complex neighborhood $\cO$ of $L=(-1,1)$ and satisfying $\phi_\tau(0)=1$, and a number $\lambda$, such that the following holds: 

\medskip

\begin{itemize}

\item[1)] $\phi_\tau$ solves the equation $(\ref{tau_equation})$ on ${\cO}$;

\medskip

\item[2)] $\phi_\tau$ has a unique quadratic critical point on $\cO$ at some real: $\phi_\tau(c+z)=O(z^2)$;

\medskip

\item[3)] the two {inverse branches} $\eta$ and $\zeta$ of $\phi_\tau$ {can be factorized} as
\bigskip
\begin{equation}\label{inverse_branches}
\eta(z)=u(T(-\sqrt{L(z)})), \quad \zeta(z)=u(T(\sqrt{L(z)})),
\end{equation}
$T$ and $L$ are  affine, and 
$$ u \in {\cA}(\cD,\cE;{\bf \mc}), \quad {\bf \mc}=\left(-{1 \over 2},0,0,1\right);$$

\medskip

\item[4)] $-0.2626 < \lambda < -0.2426$;

\medskip


\item[5)] if $\tau \equiv 0$, then $\lambda < c <0$, and there are points 
$$0<x_+<\phi_0(\lambda) <1, \quad {x_+ \over \lambda} <x_-<-1, \quad {\rm and} \quad x_*,$$
such that $\phi_0(x_\pm)=0$ and $\phi_0(x_*)=x_*$;
\item[6)] additionally, 
$$b \equiv \phi_0(c) <1+\lambda.$$
 
\end{itemize}
\end{clm}

\medskip

In this paper we will present a ``proof'' of this claim which is non-rigorous in two respects. 

First, it involves uniformization  (construction of a conformal isomorphism to the unit  circle) of the sets $\cD$ and $\cE$. This uniformization is straightforward for $\cD$, but can not be put in a simple form in the case of $\cE$. In this paper we will approximate the restriction  of the uniformizing coordinate of $\cE$ to a set compactly contained in $D(I_2,\theta_2)$ by that of $D(I_2,\theta_2)$. A rigorous construction of the uniformizing coordinate of $\cE$ would  require implementing computer bounds on a numerical approximation using, for instance the Schwarz-Christoffel formula, which is an undertaking in itself.

Second, the ``proof'' relies on verification of several inequalities of the form $f({\bf x})>0$ where $f$ is some explicit function and ${\bf x}$ ranges over a set in $\fR^d$ ($d=4$ or $5$). At this time this verification is done only numerically, its interval arithmetic implementation will be a subject of a future work.

We have also implemented the renormalization operator $\cRG$ as well as $\cK_\eps$ numerically, and have found good approximations of their fixed points by iterating Newton maps for these operators. The following summarizes our findings.

\medskip

\begin{mainobs}\label{main_obs}

\begin{itemize}
\item[]

\item[1)] The following is an approximate Taylor series for the $\phi_0$:

\begin{eqnarray}
\nonumber \phi_0(y) &\approx&  1.0 -1.9719 \times 10^{-1} y\phantom{^2} -9.2103 \times 10^{-1} y^2 -3.1550 \times 10^{-2} y^3 \\
\nonumber  &\phantom{\approx}&  \phantom{1.0}+2.5252 \times 10^{-2} y^4 -5.6774  \times 10^{-4} y^5 +1.6209 \times 10^{-5} y^6 \\
\nonumber &\phantom{\approx}&  \phantom{1.0}+5.5554 \times 10^{-6} y^7-2.5832 \times 10^{-6} y^8 +5.1783 \times 10^{-8} y^9  \ldots \\
\nonumber \lambda &\approx& -0.25014 \ldots
\end{eqnarray}


\medskip

\item[2)]
The operator ${\cK}_0$ is hyperbolic at $\phi_0$ with a local stable manifold of codimension $2$ and the two eigenvalues outside of the unit  circle given by
$$ \delta_1=8.70052...  , \quad \delta_2={1 \over {\lambda}}.$$

$\delta_2$ is the eigenvalue of the expanding eigenvector corresponding to coordinate translations.

\medskip

\item[3)]
The operator ${\cK}_\epsilon$ has a fixed point $\varphi_\eps$ for all $\epsilon$ in some neighborhood of zero.

\medskip

\item[4)]
The operator ${\cRG}$ has a fixed point $\epsilon^*$ whose approximate power series is given by
\begin{eqnarray}
\nonumber \epsilon^*=
&& -2.668 \times 10^{-2}
  +5.477 \times 10^{-2} y\phantom{^4 y^4}
  -1.385 \times 10^{-2} y^2\phantom{y^4}
  -3.543 \times 10^{-2} y^3\phantom{y^4}\\
\nonumber 
&& \phantom{-2.106 \times 10^{-2}}
  +2.187 \times 10^{-2} y^4\phantom{y^4}
  +5.504 \times 10^{-2} x\phantom{^4}\phantom{y^4}
  -6.029 \times 10^{-2} x y\phantom{^4}\phantom{^4}\\
\nonumber &&  \phantom{-2.106 \times 10^{-2}}
  +5.780 \times 10^{-3} x y^2 \phantom{^4}
  -5.777 \times 10^{-4} x y^3 \phantom{^4}
  -3.015 \times 10^{-2} x^2   \phantom{y^2}\\
\nonumber &&  \phantom{-2.109 \times 10^{-2}}
  +6.096 \times 10^{-3} x^2 y\phantom{^4}
  -9.384 \times 10^{-4} x^2 y^2
  +1.927 \times 10^{-3} x^3\phantom{y^4}\\
\nonumber  && \phantom{-2.109 \times 10^{-2}}
  -6.256 \times 10^{-4} x^3 y\phantom{^4}
  -1.444 \times 10^{-5} x^4\phantom{y^4}+\ldots
\end{eqnarray}
while
\begin{eqnarray}
\nonumber \phi_{\epsilon^*}&=&1.0-1.9468 \times 10^{-1} y\phantom{^4} -9.1492 \times 10^{-1} y^2-3.6943 \times 10^{-2} y^3\\
\nonumber &&\phantom{1.0}+2.2002 \times 10^{-2} y^4-5.9077 \times 10^{-4} y^5+1.4033 \times 10^{-5} y^6\\
\nonumber &&\phantom{1.0}+7.5534 \times 10^{-6} y^7-4.8304 \times 10^{-6} y^8+8.5963 \times 10^{-8} y^9 \ldots,\\
\nonumber \lambda^*&=&-0.24888...
\end{eqnarray}

\item[5)] The operator ${\cRG}$ is a contraction in a neighborhood of $\epsilon^*.$

\medskip

\item[6)]
The operator ${\cK}_{\epsilon^*}$ is hyperbolic at its fixed point $\varphi_\eps^*$ with the local stable manifold of codimension $2$, and the two expanding eigenvalues  given by
$$\delta_1=8.66247..., \quad \delta_2={1 \over {\lambda^*}}.$$ 
$\delta_2$ is the eigenvalue of the expanding eigenvector corresponding to coordinate translations.

\end{itemize}
\end{mainobs}

\section{Towards a complete proof of existence of a renormalization fixed point for are-preserving maps}\label{towards}

Let $\varphi_\eps$ be the solutions of $(\ref{new_eq})$, and let $s_{\phi,\epsilon}$ be a generating function of the form $(\ref{s})$ and $\tilde{s}_{\phi,\eps}$ be its renormalization. Ultimately, we would like to show that the map
$$
\eps(x,y) \mapsto \tilde{\eps}(x,y) \equiv  \tilde{s}_{\varphi_\eps,\eps}(x,y)-x +\varphi_\eps(y)
$$
is a {\it contraction} for small $\eps$.

We propose the following scheme to achieve this goal

Choose two intervals $I_x  \supset (\lambda/2,b)$ and $I_y  \supset (c,1)$, and set $\Omega \equiv I_x \times I_y$. Let $\cO(\Sigma)$, $\Sigma=I_y$ or $\Sigma=\Omega$, be the space of real-analytic functions on the open set $\Sigma$. Let $\| \cdot \|_\Sigma$ signify the sup-norm in this space.  The so called {\it Epstein} class of  functions  $\phi$ factorizable as in part $3)$ of Claim $\ref{main_thm}$ will be denoted by $\Ep(\cD,\cE,{\bf \mc})$.

\begin{itemize}
\item[1)] Let $\eps_0 \in \cO(\Omega)$ be such that $\| \eps_0 \| \le \nu_0 < {|\lambda| \delta \over 2}$. The last inequality implies that $\|\omega_{{\eps_0},\phi}\|_{I_y} \le \delta$ for all $\varphi \in \Ep(\cD,\cE,{\bf \mc})$. Since $\Ep(\cD,\cE,{\bf \mc})$ is compact, the operator
$$\varphi \mapsto \phi_{\omega_{{\eps_0},\varphi}}$$
where, as before, $\phi_\tau$ is the solution of $(\ref{tau_equation})$, has a fixed point $\varphi_{\eps_0}$ which solves $(\ref{new_eq})$ for $\eps \equiv \eps_0$. 

\item[2)] Set
$$z_0(x,y) \equiv \psi_{\eps_0} \left({\lambda \over 2 } (x+y) \right),$$
where $\psi_{\eps_0}$ be the inverse of $\varphi_{\eps_0}$ on $I_y \cap \field{R}_+$.

\item[3)] Notice, that the midpoint  equation $(\ref{new_midpoint})$  can be solved by iteration
$$z_{k+1}(x,y)=\varphi^{-1}_\eps \left({\lambda \over 2 } (x+y)+{1 \over 2} (\eps(\lambda x,z_k(x,y) ) +\eps(\lambda y,z_k(x,y)) ) \right).$$

Therefore, set 
$$z_1(x,y)=\psi_{\eps_0}\left({\lambda \over 2 } (x+y)+{1 \over 2} (\eps_0(\lambda x,z_0(x,y) ) +\eps_0(\lambda y,z_0(x,y) ) ) \right).$$

\item[4)] Set
$$\eps_1(x,y) \equiv s_{\varphi_{\eps_0} ,\eps_0} ( z_1(x,y),\lambda y )-x +\varphi_{\eps_0}(y)$$

\item[5)] Generally, at k-th step, set 
\begin{eqnarray}
\nonumber z_{k+1}(x,y) &\equiv & \psi_{\eps_k}\left({\lambda \over 2 } (x+y)+{1 \over 2} (\eps_k(\lambda x,z_k(x,y) ) +\eps_k(\lambda y,z_k(x,y) ) ) \right),\\
\nonumber \eps_{k+1}(x,y) &\equiv& s_{\varphi_{\eps_k} ,\eps_k} ( z_{k+1}(x,y),\lambda y )-x +\varphi_{\eps_k}(y)
\end{eqnarray}

\item[6)] Use {\it a-priori} bounds on $\psi_{\eps_k}$ and its derivative to estimate 
$$\| z_{k+1}-z_k\|_\Omega \le \eta_k, \quad \|\varphi_{\eps_{k+1}}-\varphi_{\eps_k} \|_{I_y} \le \upsilon_k, \quad \|\eps_{k+1} -\eps_k \|_\Omega \le \nu_k.$$

Show that  $\sum_{k=0}^{\infty} \eta_k \le \infty$,  $\sum_{k=0}^{\infty} \upsilon_{k+1} \le \infty$, and, most importantly, 
$$\sum_{k=0}^{\infty} \nu_k \le  {|\lambda| \delta \over 2},$$
which implies that the solution $\varphi_{\eps_k}$ exists and  is in $\Ep(\cD,\cE,{\bf \mc})$ for all $k$. 

\end{itemize}

We should mention that an essential ingredient for step $6)$ of the scheme is continuity of $\varphi_\eps$ in $\eps$. As it will be apparent from the discussion in the following Sections, obtaining an estimate of the type
$$\| \varphi_\eps-\varphi_{\eps'} \|_{I_y}  \le C \|\eps-\eps' \|_\Omega$$
is rather tricky, and amounts to computing a similar bound for the diffeomorphic part $u$ of the factorized inverse, as well as $\lambda$.  

We are currently working on the details of this scheme. It should be noted that the scheme is clearly convergent numerically for initial $\eps_0$ with $\| \eps_0 \|_\Omega$ much larger than  $|\lambda| \delta / 2$.

\section{Factorization of inverse branches} \label{reduction}

We will now start preparing a demonstration of Claim $\ref{main_thm}$.

We will look for the solution of $(\ref{tau_equation})$ within a class of functions which are unimodal on some interval $I  \equiv (a,d) \ni \{0,1\}$, that is they have a unique critical point on $I$, and that this critical point $c$ is quadratic in the sense that 
$$\phi_{\tau}(y)=O((y-c)^2).$$

We will now proceed to derive equations that the two inverse branches of such $\phi_{\tau}$ would satisfy. 

Write 
$$\phi_{\tau}(y)=b-g(y-c), \quad b \equiv \phi_{\tau}(c),$$
then $(\ref{tau_equation})$ can be written as 
\begin{equation}\label{g_eq_2}
 g=F \circ g \circ \xi + id-\tau \circ (id +c),
\end{equation}
where
\begin{eqnarray}
\nonumber F(y)&=&b+c-{2\over \lambda} (b-g(b+c-y)), \\
 \nonumber \xi(y)&=&\lambda y +c(\lambda-1).
\end{eqnarray} 

Denote $h$ and $f$ the two inverse branches of $g$:
$$h: (0,g(d-c)) \mapsto (0,d-c), \quad f: (0,g(a-c)) \mapsto (a-c,0).$$

The ``inverse'' of  $(\ref{g_eq_2})$ is the following set of equations for the inverse branches:
\begin{eqnarray}
\label{branch_1} f \circ F^{-1} \circ (id -h+\tau  \circ (h+c))&=&\xi \circ h, \quad {\rm on} \quad (E,g(d-c)) \\
\label{branch_2} h \circ F^{-1} \circ (id -h+\tau  \circ (h+c))&=&\xi \circ h, \quad {\rm on} \quad (0,E),\\
\label{branch_3} h \circ F^{-1} \circ (id - f + \tau \circ (f+c))&=&\xi \circ f,  \quad {\rm on} (0,g(a-c)),
\end{eqnarray}
where
$$E\equiv e^2, \quad e \equiv -\sqrt{g\left({c \over \lambda}-c\right)}.$$ 


\begin{figure}
 \begin{center}
\begin{tabular}{c c c}
$\!\!\!\!\!\!\!$ \resizebox{50mm}{!}{\includegraphics{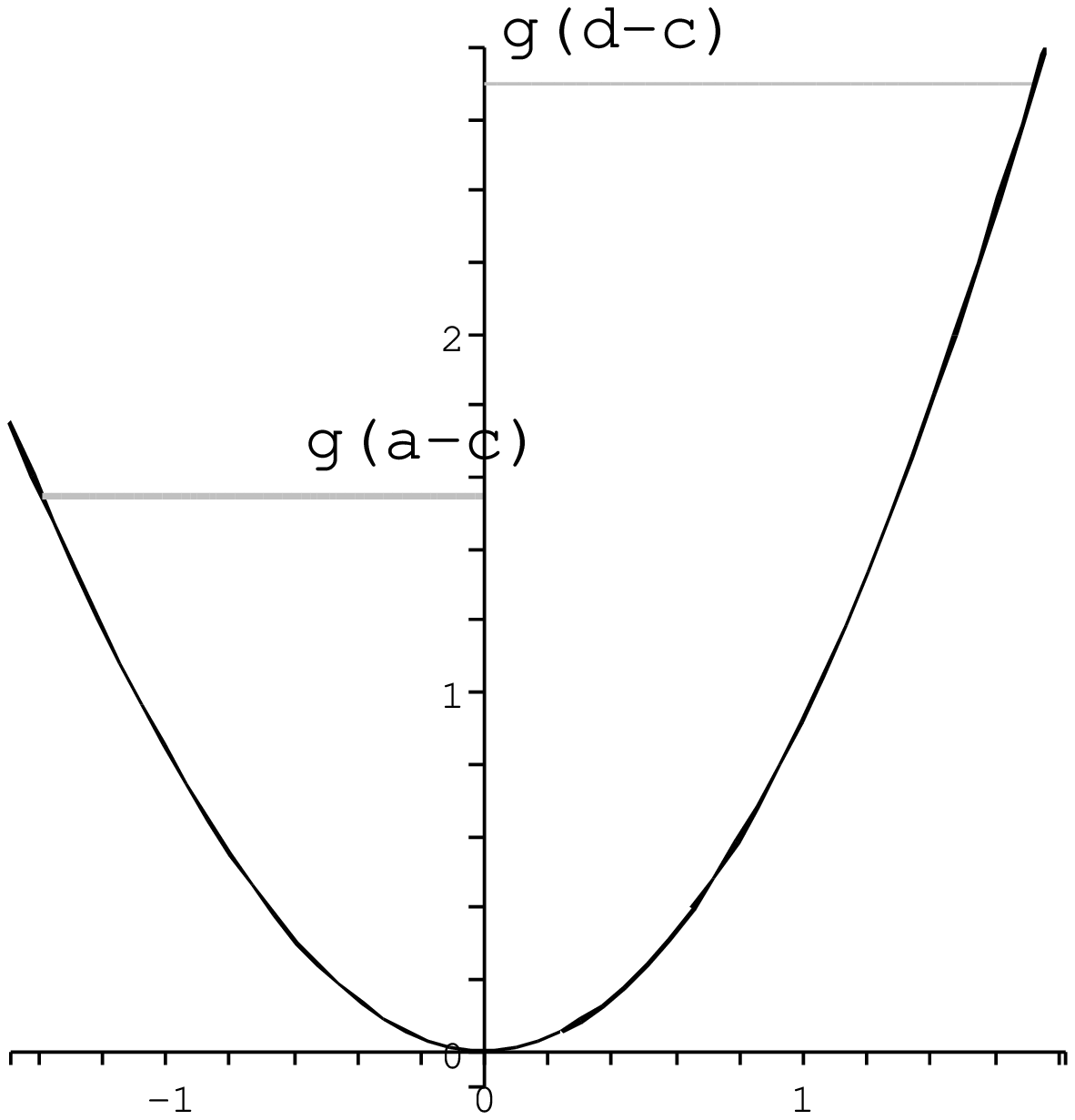}} &  $\!\!\!\!\!$ \resizebox{50mm}{!}{\includegraphics{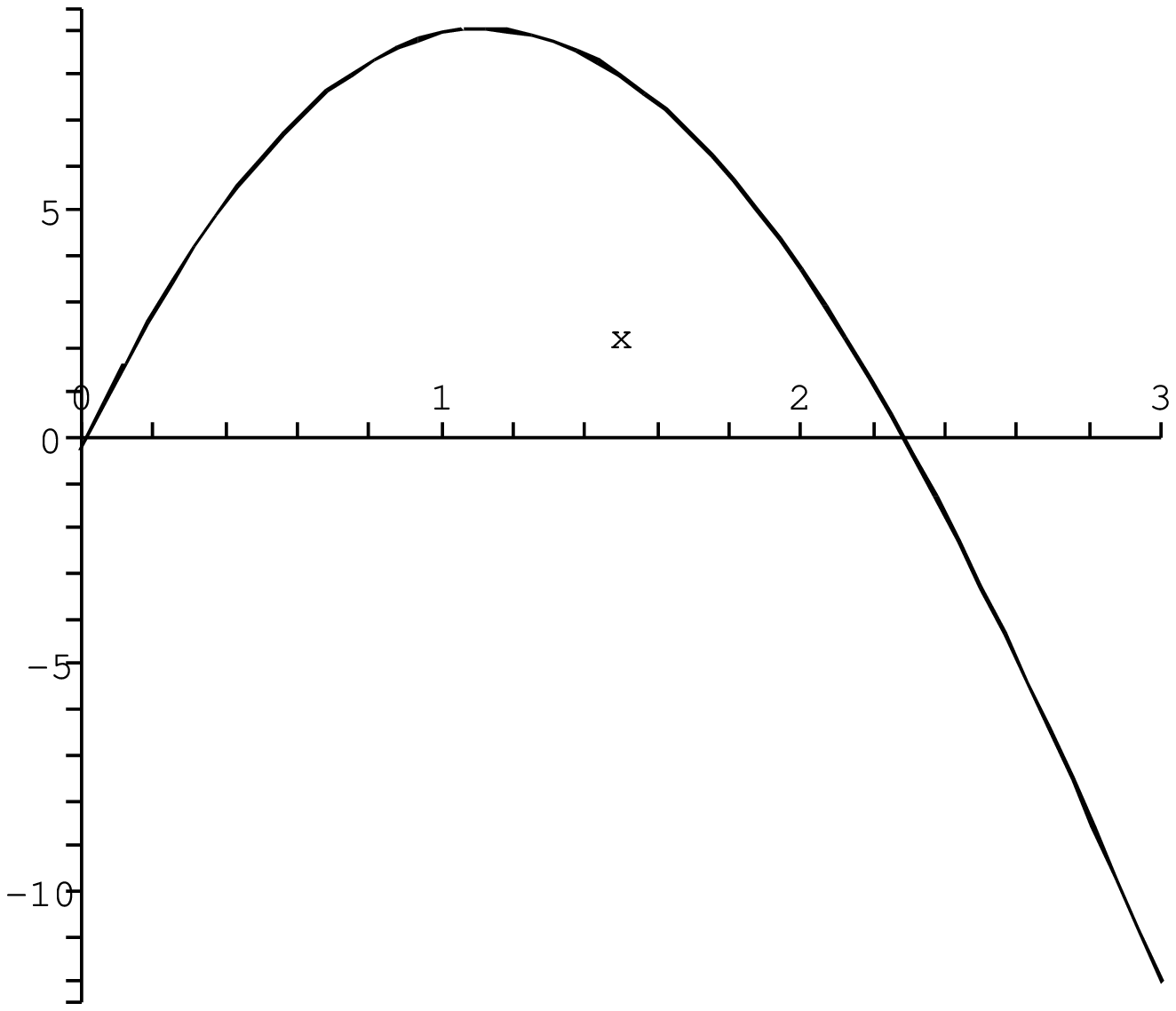}} & $\!\!\!\!\!$ \resizebox{50mm}{!}{\includegraphics{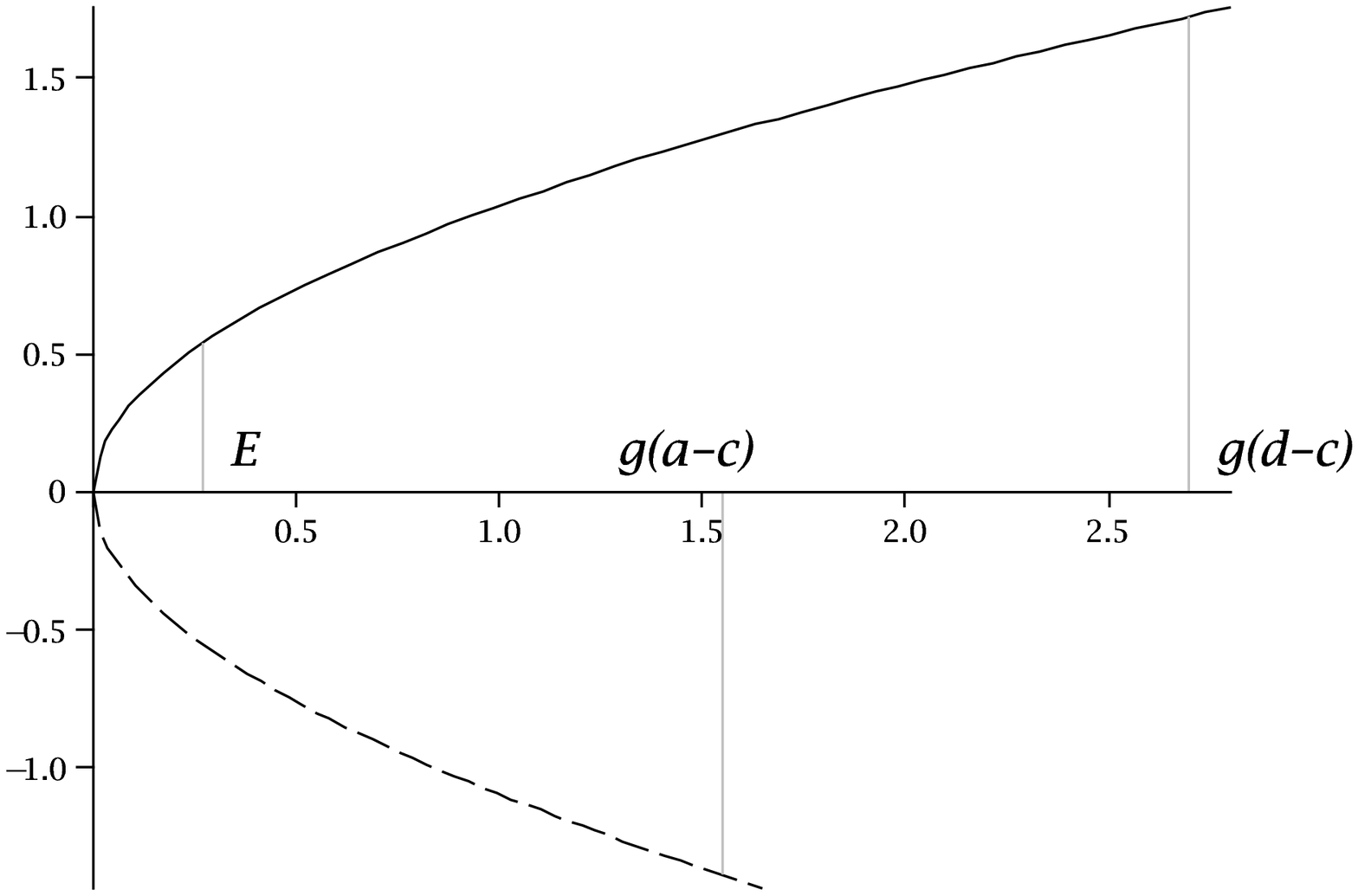}}\\
$\!\!\!\!\!\!\!$ a) &$\!\!\!\!\!$ b) & $\!\!\!\!\!$ c)

\end{tabular}
\caption{Function $g$ (in a)), function $F$ (in b)) and inverse branches $h$ and $f$ (in c)) for the solution $\phi$ of the equation $\phi(x)=2 \lambda^{-1}  \phi(\phi(\lambda x))-x$.} 
\end{center}
\end{figure}

We will  look for solutions of $(\ref{tau_equation})$ within the {\it Epstein class}:
\begin{equation}\label{L_Epstein_class}
\phi(y)=U(y)^2,
\end{equation}
where $U$ is  diffeomorphism of $I$ onto its image, and we will write
\begin{equation}\label{factorization}
{h=v \circ - \circ  s, \quad f=v \circ s}, 
\end{equation}
where $v$ is a diffeomorphism on $K \equiv (-\sqrt{g(d-c)},\sqrt{g(a-c)})$, $s(x) \equiv \sqrt{x}$ (the principle square root) and $-(x) \equiv -x$. A similar factorization has been used in \cite{Sul} and \cite{LY}, however unlike the authors of those works we will not make any assumptions on the univalence of $v$ in some neighborhood of $K$; below we will choose a space of $v$'s as functions holomorphic on a neighborhood of $K$, and we will point out  a specific obstruction to univalence. 

With this factorization equations $(\ref{branch_1})$, $(\ref{branch_2})$ and $(\ref{branch_3})$ become
\begin{eqnarray}
\label{linearizer}\xi \circ v &=& v \circ V, \\
V(x)&=&-\left[F^{-1} (x^2-v(x)+\tau(v(x)+c)  )\right]^{1 \over 2}, \quad  x \in [e,\sqrt{g(a-c)}), \\
V(x)&=&\left[F^{-1}  (x^2- v(x)+\tau(v(x)+c))\right]^{1 \over 2}, \quad  x  \in (-\sqrt{g(d-c)},e).
\end{eqnarray}

\begin{figure}
 \begin{center}
\begin{tabular}{c c c}
 \resizebox{45mm}{!}{\includegraphics{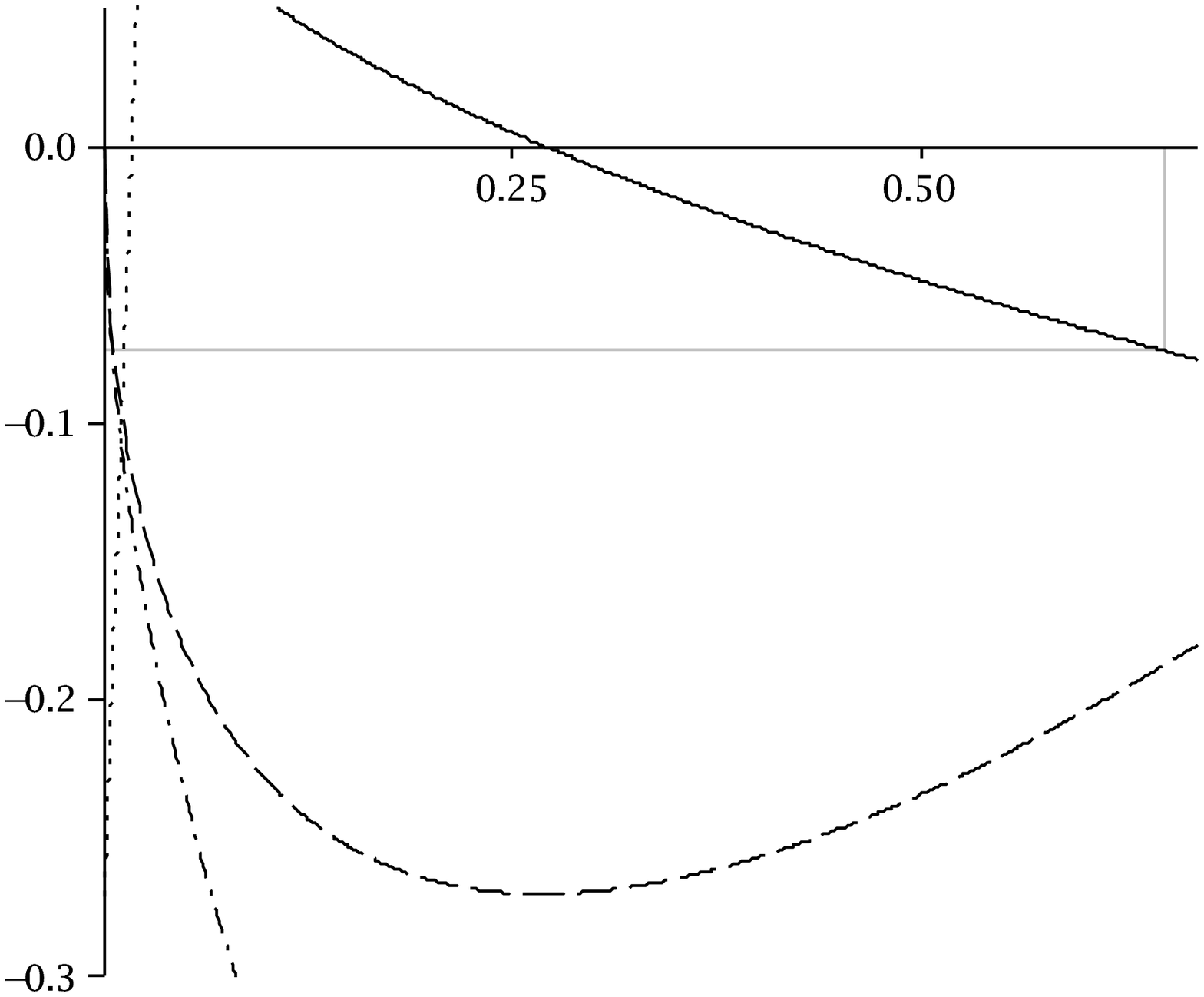}} &   \resizebox{45mm}{!}{\includegraphics{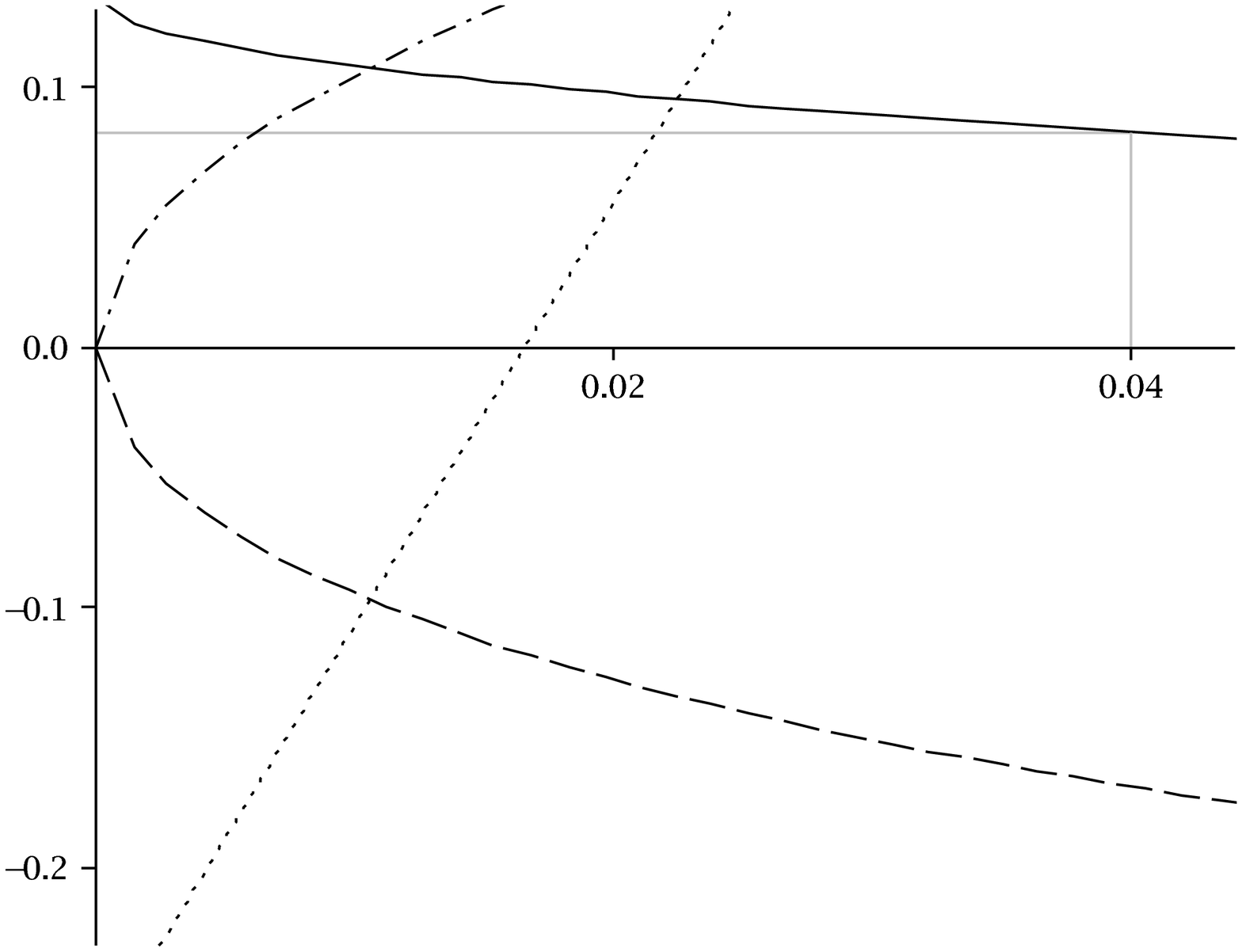}} & \resizebox{50mm}{!}{\includegraphics{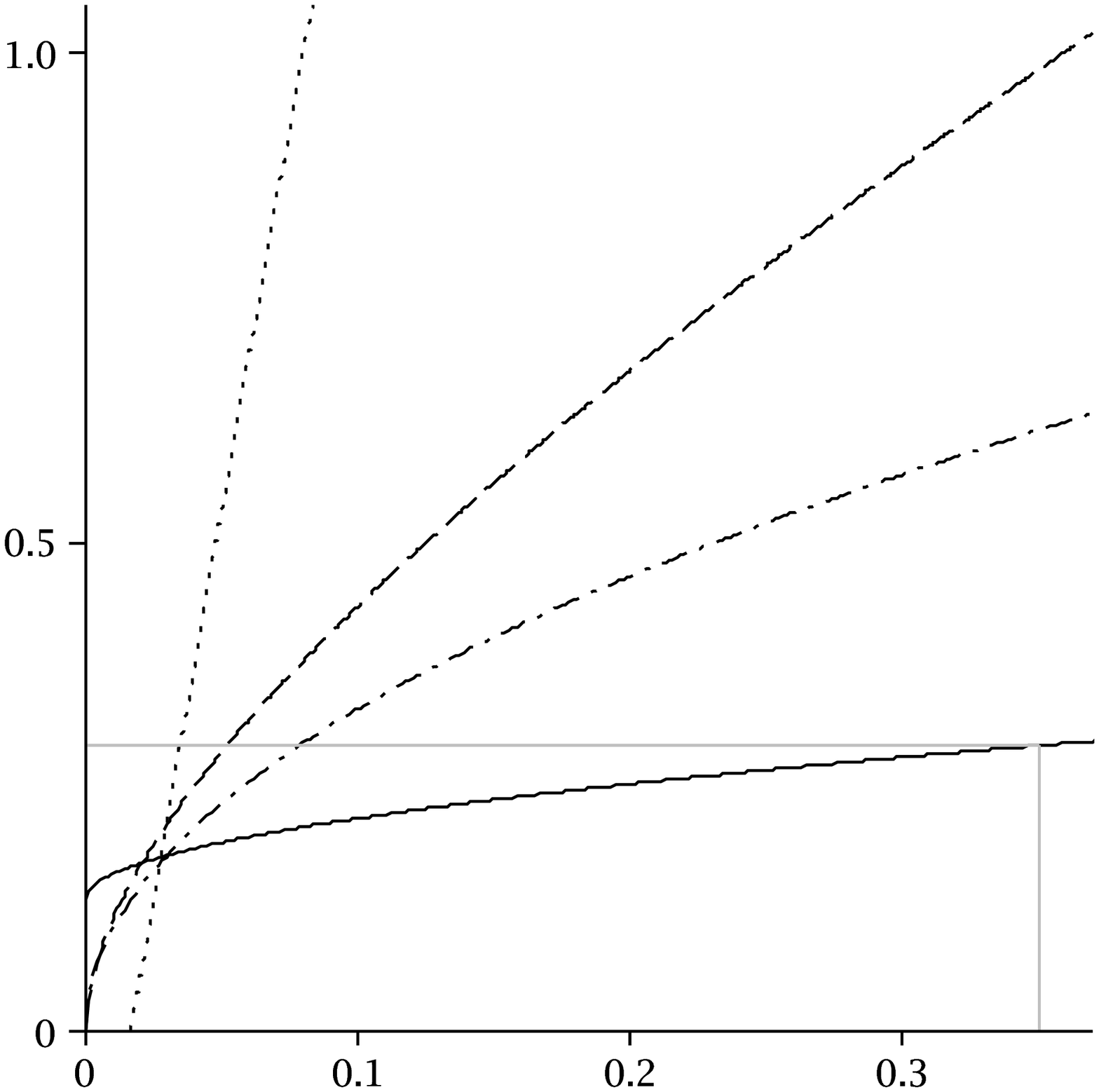}} \\
a) & c) & e) \\
 \resizebox{45mm}{!}{\includegraphics{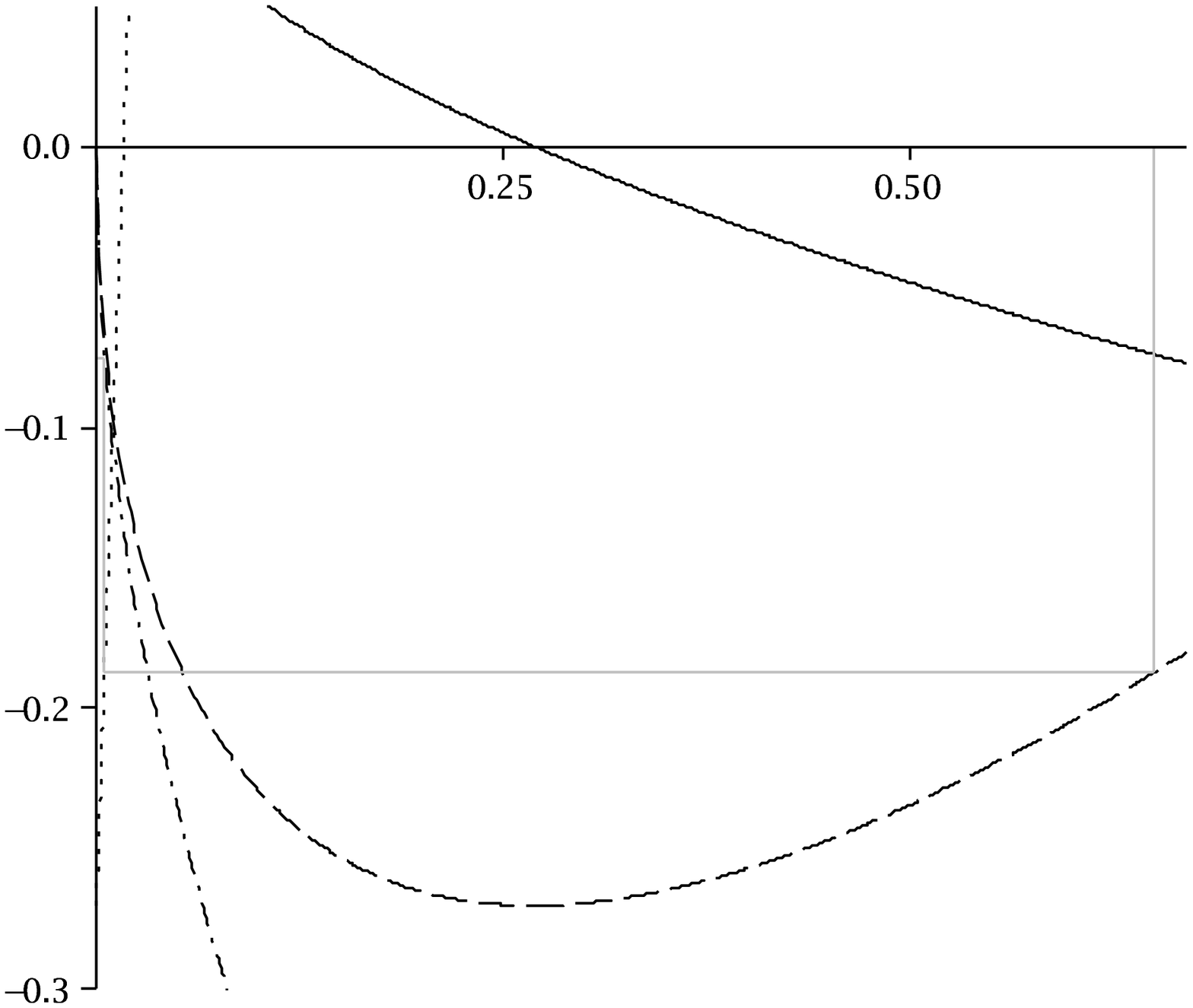}} &  \resizebox{45mm}{!}{\includegraphics{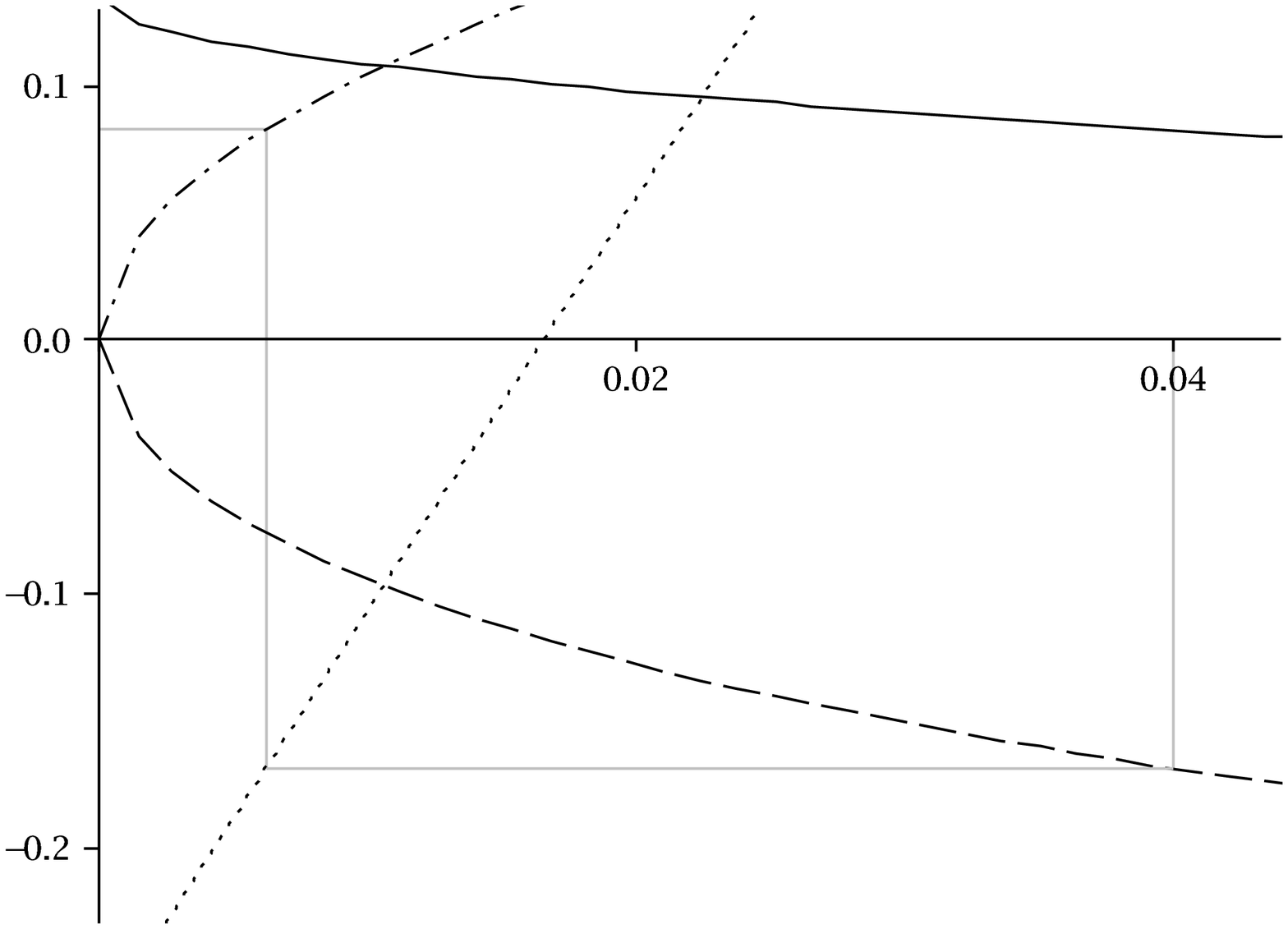}} & \resizebox{45mm}{!}{\includegraphics{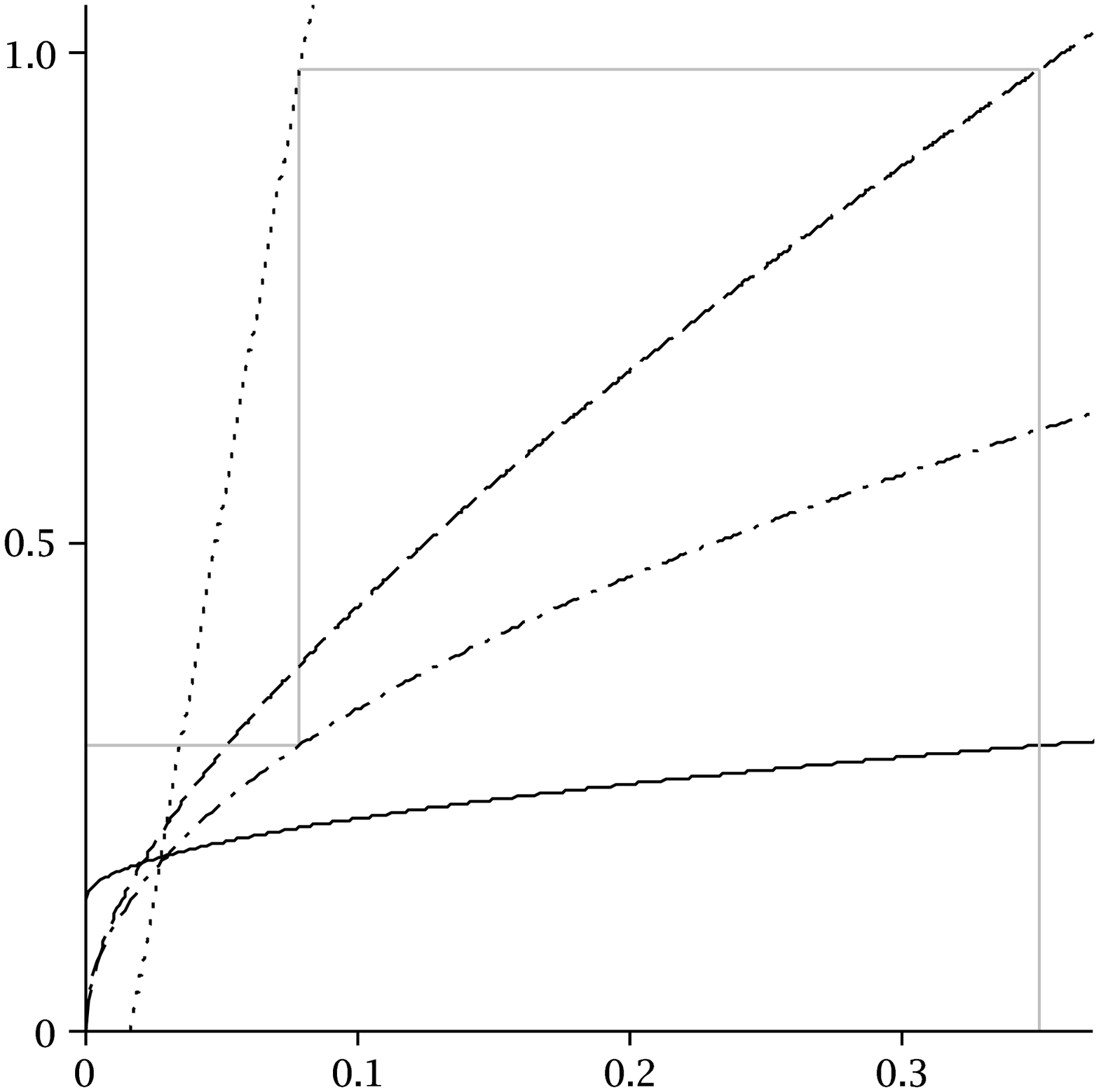}}\\
b) & d) & f)
\end{tabular}
\caption{Combinatorics in the equalities $(\ref{branch_1})$ -$(\ref{branch_3})$. {\bf Equality $\mathbf{(\ref{branch_1})}$, a)  and  b):} function $\xi \circ h$ is given in a solid line, $id-h$ -- in dash, $F$ -- in dot, $f$ -- in dash-dot; the image of the point under the right hand side of the equality is shown in $a)$, under the left hand side -- in $b)$.  {\bf Equality $\mathbf{(\ref{branch_2})}$, c)  and  d):} function $\xi \circ h$ is given in solid, $id-h$ -- in dash, $F$ -- in dot, $h$ -- in dash-dot; the image of the point under the right hand side of the equality is shown in $c)$, under the left hand side -- in $d)$.  {\bf Equality $\mathbf{(\ref{branch_3})}$, e)  and  f):} function $\xi \circ f$ is given in solid, $id-f$ -- in dash, $F$ -- in dot, $h$ -- in dash-dot; the image of the point under the right hand side of the equality is shown in $e)$, under the left hand side -- in $f)$.}
\end{center}
\end{figure}

Below, we will choose a compact functional space for $v's$, and we will demonstrate that the operator $v \mapsto \xi^{-1} v \circ  V$ is defined and is continuous on this functional space; existence of a fixed point $v^*$ (and hence, existence of $\phi_{\tau}$) will follow from the Schauder-Tikhonov Theorem which guarantees existence of a fixed point for a continuous operator on a compact set.

\section{An operator on a compact space}

We will now formally introduce an operator which will be later shown to be defined on ${\cA}(\cD,\cE,{\bf \mc})$ for some choice of $\cD$, $\cE$ and 
$${\bf \mc}=\left(-{1 \over 2}, 0, 0,1 \right).$$

The operator is defined through the following sequence of steps.

\begin{itemize}
 \item[i)] Given $ u \in {\cA}(\cD,\cE,{\bf \mc})$, and a function $\tau$, holomorphic on $\cE \ni 0$, real-valued on $\fR$ and satisfying $\tau(0)=0$, find $b, \lambda$  and $e$ from the following set of equations:
\begin{eqnarray}
\label{e_equation} -2 e  &=&  \alpha(b,\lambda) (1-\tau'(u(\Tble(e)) )) u'(\Tble(e)), \\
\label{l_equation} \lambda   &=&   u \left(\Tble \left(  \left[ { b - u\left( \Tble  \left[  -\sqrt {b- {  \lambda \over 2}-{\lambda^2 \over 4}-\tau(1)} \right] \right) } \right]^{1 \over 2} \right)\! \right),\\
\label{b_equation} b  &=&  u \left( \Tble \left[  \sqrt{b- { \lambda \over 2} ( b-e^2+u(\Tble (e)) -\tau(u(\Tble(e)))} \right] \right),
\end{eqnarray}
where $\alpha$,  and additional functions $\beta$ and $\gamma$ are given by
\begin{equation}
\nonumber \alpha(b,\lambda)= { 1 \over 2 \beta(b,\lambda)-2 \gamma(b)}, \quad \beta(b,\lambda) =\sqrt{ b-{\lambda \over 2}}, \quad  \gamma(b)=\sqrt{b-1},
\end{equation}
and
\begin{equation}
\nonumber \Tble(x)= -\alpha(b,\lambda) (x+\beta(b,\lambda)).
\end{equation}

\item[ii)] Define for all $x \in T^{-1}_{b,\lambda}(\cD \cap \fR)$ 
\begin{equation}\label{V_equation}
 \Veut(x) = {\rm sign}(e-x) \left[{b-u  \left( \Tble  \left( - \left[ {w(\Tble(x))} \right]^{1 \over 2}\right)\right)} \right]^{1 \over 2},
\end{equation}
where
\begin{equation}\label{w_function}
 w (z) ={b- {\lambda \over 2} \left(b-\Tble^{-1}(x)^2+  u(x) +\tau(u(x))\right)}
\end{equation}

We will demonstrate that there is a choice of $\cD$ and $\cE$ such that  $\Veut$ extends to a holomorphic function on $T^{-1}_{b,\lambda}(\cD)$.
 
\item[iii)] Set

\begin{equation} \label{T_op}
{\cT}_{\tau}[u](\Tble(z)) \equiv \lambda^{-1} u(\Tble(\Veut(z))) .
\end{equation}

The operator $\cT_{0}$ will be also denoted by $\cT$.

\end{itemize}

\begin{rmk}
{\rm Notice, that $\gamma=-\sqrt{b-1} \in (e,0)$ is the fixed point of $\Veut$.}
\end{rmk}

\begin{rmk}
{\rm  The normalization conditions $(\ref{e_equation})$--$(\ref{b_equation})$ ensure that $\Veut$ is differentiable at $e$, and that } 
\begin{equation}
\nonumber \oT[u](-1/2)=0, \quad \oT[u](0)=1.
\end{equation}
\end{rmk}

\begin{rmk}\label{fact_inverse}
{\rm The  function $u$ is related to functions $v$, $\psi$, $h$ and $f$ appearing in Section $\ref{reduction}$ through the following equations:
\begin{eqnarray}
  \nonumber v(x)&=&u(-\alpha(x+\beta))-c, \\
  \nonumber h(x) &\equiv & =\psi(b-x)-c=u(\alpha (\sqrt{x}-\beta))-c, \quad x \in \left(0, \left[T^{-1}_{b,\lambda}(r)\right]^2\right), \\
  \nonumber f(x) &\equiv & u(\alpha (-\sqrt{x}-\beta))-c, \quad x \in  \left(0, \left[T^{-1}_{b,\lambda}(l)\right]^2\right),
\end{eqnarray}
}
\end{rmk}

We will show that for small $\tau$, there is a choice of $\cD$ and $\cE$ such that that ${\cT}_{\tau}[u]  \in \cA(\cD,\cE,{\bf \mc}  )$ whenever $u \in  \cA(\cD,\cE,{\bf \mc}  )$. By compactness of the set $\cA(\cD,\cE,{\bf \mc}  )$ there is a function $u^*_{\tau} \in \cA(\cD,\cE,{\bf \mc}  )$ such that ${\cT}_{\tau}[u^*_{\tau}]=u^*_{\tau}$, which is equivalent to the set of equations $(\ref{branch_1})-(\ref{branch_3})$. In particular, $u^*_{\tau}$ is the ``factorized inverse'' (in the sense of Remark \ref{fact_inverse}) of a solution  of the equation $(\ref{tau_equation})$.

\begin{rmk}
{\rm 
Before we proceed with the proofs, we would like to emphasize two crucial difficulties that have forced us to modify the standard techniques that are commonly used to control inverse branches of unimodal maps (cf. \cite{Eps1}, \cite{Eps2}, \cite{Sul}, \cite{LY}).  

\begin{itemize}
\item[i)] The terms $y$ and $\tau(y)$ in the equation $(\ref{tau_equation})$ are responsible for the appearance of the terms $u(\Tble(x))$ and $\tau(u(\Tble(x))))$ in $(\ref{V_equation})$. The effect of these terms is that one looses the benefit of needing to estimate $u$, every time it enters the expression for $V_{\lambda,\tau}$, only on a compact subset of its domain where one can use {\it a-priori} bounds. These terms do not appear in the Feigenbaum case  where this difficulty is absent. {\it In the current situation one can not but  make assumptions on the range of $u$, and show that these assumptions are reproduced.}

 \item[ii)] Another effect of terms $u(\Tble(x))$ and $\tau(u(\Tble(x))))$ in $(\ref{V_equation})$ is that derivative
$$\oT[u]'(z)=-\lambda^{-1} u'(\Tble(\Veut(\Tble^{-1}(z))) \alpha \Veut(\Tble^{-1}(z))'$$
can become zero since
$$\Veut(\Tble(z))'=\ldots \times {1 \over \Veut(\Tble^{-1}(z))} \left({2 \over \alpha}  \Tble^{-1}(z)+u'(z)-\tau'(u(z)) u'(z)\right)$$
can be zero. Notice, that $\Veut(\Tble(z))'$ is not zero at $e$ where an application of the L'Hopital's rule shows that the derivative is finite. However, it can be zero at other points on the real line where $2\alpha^{-1} \Tble^{-1}(z)+u'(z)-\tau'(u(z)) u'(z)$ is zero. This would totally destroy the argument since a function $\tilde{u} \equiv \oT[u]$ whose derivative is zero somewhere in the real slice of its domain generally is not in $\cA(\cD,\cE,{\bf \mc}  )$, in particular $\tilde{u}(\cD \cap \fC_\pm) \nsubseteq \overline{\cE \cap \fC_\pm})$. 

{\it We will deal with this problem by assuming an upper bound on the derivative $u'$ in the ``problematic'' subinterval of the real slice of $\cD$ so that $2 \alpha^{-1} \Tble^{-1}(z)+u'(z)-\tau'(u(z)) u'(z)$ is guaranteed to be nonzero, and we will demonstrate that this bound is reproduced.}

\end{itemize}
 
}
\end{rmk}
\bigskip

\section{Detailed statement of the Main Claim}\label{epsneq0}

Let $S_{\vartheta}$ denote the sector in $\field{C}$ of angle $2 \vartheta$ with the vertex at point $-1$, containing $(-\infty,-1)$ and symmetric with respect to the real axis. We denote

$$\field{C}_{-1,\vartheta}=\field{C}_1 \setminus S_\vartheta.$$

In what follows, we will make the following choices: 
$$\cD=D(I_1,\theta_1), \quad \cE=\cE_\vartheta \cup D(I_3,\theta_3) \cup D(I_4,\theta_4),$$
where
$$ \quad \cE_\vartheta=\tilde{\Theta}_2( \field{C}_{-1,\vartheta}) \subset  D(I_2,\theta_2),$$
and $I_k=(l_k,r_k)$ and $\theta_k$ are as in $(\ref{I1})$---$(\ref{I3})$, $\vartheta=0.24$ and $\tilde{\Theta}_2$ is a normalized conformal isomorphism of $\cD=D(I_2,\theta_2)$ and $\field{C}_1$ (see below).  We will consider the corresponding space $\cA(\cD,\cE,{\bf \mc})$.

The point of considering such rather peculiar subset of a Poincar\'e neighbourhood as the target  set is that the factorized inverse $u$ of the solution $\phi_\tau$ of $(\ref{tau_equation})$ generally maps a symmetric domain into a non-symmetric one. In our numerical experiments  we have determined that the above choice $\cE$ seems to be a relatively good approximation of $u(\cD)$, being, at the same time, the set for which one can give a relatively fair  approximation of the conformal isomorphism with $\field{C}_1$.

Denote 
$$\cD_\pm \equiv \cD \cap \field{C}_\pm,$$
and let 
$$\theta' \mapsto \partial \cD_+(\theta'), \quad \theta' \in \left( \theta_1-{\pi \over 2},{3  \over 2} \pi-\theta_1 \right),$$
be the standard parametrization of the circle arc $\partial D_+(I_1,\theta_1)$. In this paper we will use an affine rescaling of $\theta'$ to the interval $(0,\pi)$:
$$\theta \mapsto \partial \cD_+(\theta), \quad \theta={1 \over 2 \pi - 2 \theta_1  } \left( \theta'-\theta_1+{\pi \over 2 }\right)$$
as  our parametrization of the boundary of $\cD_+$. In a similar way, $\cD_-$ is parametrized by $\theta \in (-\pi,0)$.
  
\begin{figure}
\vspace{-2mm}
 \begin{center}
\resizebox{50mm}{!}{\includegraphics[angle=-90]{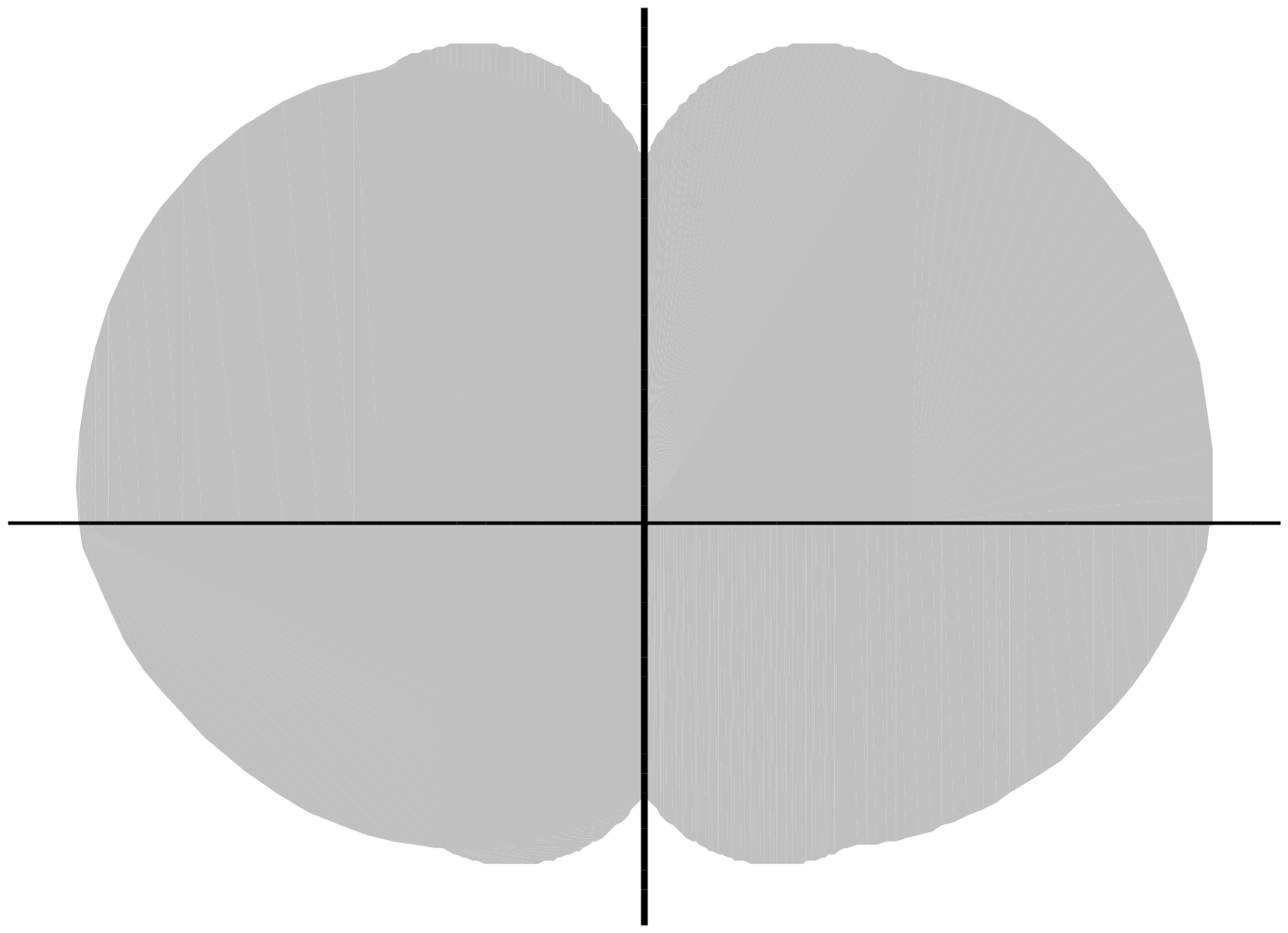}}
\vspace{-2mm}
\caption{\it Set $\cE$.}
\end{center}
\vspace{-5mm}
\end{figure}

The double slit plane $\fC_1$ is isomorphic to Poincar\'e neighbourhoods $D(I_k,\theta_k)$ via  conformal isomorphisms
\begin{equation}\label{theta_k}
\tilde{\Theta}_k \equiv q_k \circ \sigma_k \circ m_k \circ  \zeta,
\end{equation}
where
\begin{eqnarray}
\nonumber \zeta (z)& \equiv & { \sqrt{1+z}-\sqrt{1-z} \over \sqrt{1+z} + \sqrt{1-z} }, \\
\nonumber  m_k (z)& \equiv & { z+a_k \over 1 + a_k z}, \\
\nonumber \sigma_k (z)& \equiv& { (c_k+1)(1+z)^\kappa-(c_k-1) i^\kappa(1-z)^\kappa  \over (c_k+1)(1+z)^\kappa + (c_k-1) i^\kappa(1-z)^\kappa }, \\
\nonumber q_k (z)&\equiv& {r_k-l_k \over 2} z +{r_k+l_k \over 2},
\end{eqnarray}
$l_k$ and $r_k$ are the left and the right end points of intervals $I_k$,  and 
$$c_k \equiv  \Im \left( {e^{i \kappa_k \pi }+1+2 {\rm signum}(1-\kappa_k) \sqrt{e^{i \kappa_k \pi }  } \over e^{i \kappa_k \pi }-1}  \right), \quad \kappa_k \equiv 2- {2 \theta_k \over \pi}.$$

With a little bit of  work, one can check that the transformation $\zeta$ maps $\fC_1$ onto the unit disk, $m_k$ is the normalizing Moebius transformation, $\sigma_k$ maps the unit disk onto $D((-1,1),\theta_k)$, and, finally, $q_k$ maps  $D((-1,1),\theta_k)$ onto  $D(I_k,\theta_k)$.

Furthermore, the transformation

\begin{equation}
\Pi_2(z) \equiv  -2^{\vartheta \over \pi-\vartheta} (1-z)^{\pi \over \pi-\vartheta }+1
\end{equation}
maps $\fC_{-1,\vartheta}$ conformally onto $\fC_1$. Finally, we set

\begin{equation}
\Theta_1 \equiv \tilde{\Theta}_1, \quad \Theta_2 \equiv \Pi_2 \circ \tilde{\Theta_2}, \quad \Phi_1 \equiv \Theta_1^{-1}, \quad \Phi_2 \equiv \Theta_2^{-1}
\end{equation}

The constants $a_k$ in the normalizing Moebius transformations $m_k$ are defined through the conditions 
$$\Theta_1(0) =-1/2, \quad \Theta_2(0)=0.$$


A function  $u$ in $\cA(\cD,\cE,{\bf \mc})$ can be now  factorized  approximately as 
$$u= \Theta_2 \circ  f \circ \Phi_1$$
where $f \in \cA_1({\bf c})$. We emphasize that this is only approximate factorization since 
$$\cE \ne \cE_\vartheta.$$

According to Schwarz Lemma $\ref{Epstein_lemma}$, if $f \in \cA_1({\bf c})$ and an interval $J$ are such that $g(J) \subset J'$ then 
$$u( \Theta_1(D(J,\theta))) \subset \Theta_2(D(J',\theta)).$$

Furthermore, one can use the fact that $\Theta_k \arrowvert_\fR$ are monotone functions to transfer the improved Herglotz bounds $(\ref{f_bounds})$ from $\cA_1({\bf c})$ to $\cA(\cD,\cE,{\bf \mc})$:
\begin{eqnarray}
\label{U_bounds}\Um(x;t,s) &\equiv& \Theta_2\left(\mF\left(\Phi_1(x);t,s\right)\right), \\
\label{u_bounds}\um(x;t,s) &\equiv&  \Theta_2\left(\mf\left(\Phi_1(x);t,s\right)\right).
\end{eqnarray}

We have implemented bounds $(\ref{U_bounds})$--$(\ref{u_bounds})$ on the computer, and used them for a numerical verification of various conditions in our proofs.

We shall now proceed to describe a set $\cS$ of realizable derivatives $(u'\left(-{1 \over 2} \right),u'(0))$: 

\begin{lem}\label{lemma_S}
There is a convex open set ${\cS} \subset \fR^2$ , such that
$$\left(u'\left(-{1 \over 2} \right),u'(0)\right) \in \cS,$$
whenever  $u \in {\cA}(\cD,\cE,{\bf \mc})$.
\end{lem}

\noindent {\it Proof.} See Subsection $\ref{ts}$ of the Appendix for the proof. The set $\cS$ is contained in the square with sides $t_-=1.9775$, $t_+=2.0229$ and $s_-=2.011$, $s_+=2.04621$.

$\Box$

The following result shows that a set of  functions  $u \in {\cA}(\cD,\cE,{\bf \mc})$ that satisfy a certain set of conditions on the distortion of $\partial \cD$ is invariant under $\oT$.

\begin{clm}\label{central_prop}\footnote{The proof of this claim relies on several verifications of inclusion of a domain in  $\field{C}$ that continuously depends on four real parameters in another domain in $\field{C}$, and verification of several functional inequalities. These containments and inequalities have been verified numerically on a computer, however presently we have not made this computer verification rigorous (that is the programs do not use interval arithmetics yet).}


Suppose that $u  \in {\cA}(\cD,\cE,{\bf \mc})$ satisfies the following two sets of conditions.

\begin{itemize}

\item[1)] The distortion of the boundary of the domain by $u$ is bounded:

\begin{equation}
\label{cond1} u(\partial \cD(\theta) ) \in \left( \cE  \setminus D(I_5,\theta_5) \right)  \cap \{z  \!\in \! \field{C}: c_{\Im}(\theta) |\Im(\partial \cD(\theta) )|  < \Im(z) < C_{\Im}(\theta) |\Im(\partial \cD(\theta) )|  \},
\end{equation}
where
\begin{equation}
\label{cond2} I_5=\left(\max_{(t,s) \in \cS}\Um(l_1;t,s), \min_{(t,s) \in \cS}\um(r_1;t,s) \right), \quad \theta_5=1.3 \theta_2,  \quad \theta^*=1.0,
\end{equation}
and $c_{\Im}(\theta)$ and $C_{\Im}(\theta)$ are piecewise linear functions \footnote{Values of $c_{\Im}(\theta)$ and $C_{\Im}(\theta)$ are linear interpolations of the following values:
\medskip
\begin{center}
\begin{tabular}{|c | c| c| c| c| c| c| c| c| c| c| c|}
\hline
$\theta$& $0.0$ & $0.1$  & $0.2$ & $0.3$ & $0.4$ & $0.5$ & $0.6$ & $0.7$ & $0.8$ & $0.9$ & $1.0$ \cr
\hline
$c_{\Im}(\theta)$ & $2.5$ & $2.5$ &$2.5$ &$2.5$ &$2.5$ &$2.5$ &$2.5$ &$2.3$ & $2.3$& $2.2$ & $2.1$ \cr
\hline
$C_{\Im}(\theta)$ & $4.5$ & $4.5$ &$4.5$ &$4.5$ &$4.5$ &$4.0$ &$3.5$ &$2.8$ & $2.7$& $2.6$ & $2.5$ \cr
\hline
\end{tabular} 
\end{center}
}

\item[2)] The derivative of $u$ on a subinterval of the real slice of the domain is bounded:

\begin{equation}\label{der_bound}
u'(x) \le  \omega + \sigma x, \quad {\rm where} \quad  \omega=15, \quad \sigma=30,
\end{equation}
for all $x \in (0,r_1)$.

\end{itemize}

\medskip

Then, there are $\delta>0$ and $\kappa>0$, and four functions \footnote{Values of $a_i$, $A_i$, $b_i$ and $B_i$ are as follows:

\begin{tabular}{c c c c c c}
 $a_1=-5.35599046$, &  $a_2=0.139748610$, &  $a_3=10.2175899$, & $A_1=0.21030246$, &  $A_2=-0.607475331$,&  $A_3=0.556299141$, \cr
$b_1=2.55648917$, &  $b_2=0.607658324$, & $b_3=1.99870694$, &  $b_4=0.5118151460$, & $b_5=-9.20477427$, &  $b_6=-1.26473672$,\cr
$B_1=0.31811096$, &  $B_2=0.169200466$, & $B_3=-0.656080571$, &  $B_4=-0.222331242$, &  $B_5=0.638677617$,  & $B_6=1.10774246$.
\end{tabular} 
}
\begin{eqnarray}
\cL_-(t,s)& \equiv & a_1 t+a_2 s +a_3,\\
\cL_+(t,s)& \equiv & A_1 t+A_2 s +A_3, 
\end{eqnarray}
and 
\begin{eqnarray}
\cB_-(\lambda,t,s)&\equiv& (b_1 \lambda+b_2) t+(b_3 \lambda +b_4)s + b_5 \lambda +b_6,\\
\cB_+(\lambda,t,s)&\equiv & (B_1 \lambda+B_2) t+(B_3 \lambda +B_4)s + B_5 \lambda +B_6,
\end{eqnarray}
such that for any $\tau$, holomorphic on $\cE$, and satisfying 

\begin{equation}\label{tau_conds} 
\sup_{z \in \cE} |\tau(z)| \le \delta, \quad  \sup_{z \in \cE} |\tau'(z)| \le \kappa, \quad \tau(0)=0,
\end{equation}
the following holds:

\begin{itemize}

\item[i)] there is a continuous branch of solutions of equations  $(\ref{e_equation})$--$(\ref{l_equation})$; that is a triple $\mS \equiv (e,b,\lambda)$ that solves  $(\ref{e_equation})$--$(\ref{l_equation})$ and is such that the map $(u,\tau) \mapsto \mS$ is continuous for all $u \in  {\cA}(\cD,\cE,{\bf \mc})$ and $\tau$ as in $(\ref{tau_conds})$. Furthermore, it satisfies
\begin{eqnarray}
\label{e_bounds} -\gamma(b) \ge  &e& \ge -\beta(b,\lambda), \\
\label{l_bounds}  \cL_+(t,s) \ge  &\lambda& \ge \cL_-(t,s),\\
\label{b_bounds}  \cB_+(\lambda,t,s) \ge & b &\ge \cB_-(\lambda,t,s).
\end{eqnarray}

\phantom{aaa}

\item[ii)] $\oT[u]'$ also admits the bound $(\ref{der_bound})$;

\phantom{aaa}

\item[iii)]  the function $\Veut$ extends to a conformal function on $T^{-1}_{b,\lambda}(\cD)$ that maps $T^{-1}_{b,\lambda}(\cD) \cup \fC_\pm $ compactly into  $T^{-1}_{b,\lambda}(\cD) \cup \fC_\mp $;

\phantom{aaa}

\item[iv)] $\cT_{\tau}[u] \in \cA(\cD,\cE, {\bf \mc})$;

\phantom{aaa}

\item[v)] $\cT_{\tau}[u]$ satisfies condition $(\ref{cond1})$.

\end{itemize}
\end{clm}

A rather technical demonstration of this result is given in the Appendix.

\begin{rmk}
{\rm We do not demonstrate uniqueness of the solution $(e,b,\lambda)$, although this seems possible (with significantly more effort).}
\end{rmk}

We fix a branch of solutions $\mS$ (there is at least one such branch by Part $2)$, $i)$ of Claim $\ref{central_prop}$.  Together with the definition of the operator $\cT_{\tau}[u]$ this implies that this operator is {\it continuous} on the subset of  $\cA(\cD,\cE, {\bf \mc})$ of functions that satisfy the conditions $(\ref{cond1})$ and $(\ref{der_bound})$. By the Schauder-Tikhonov Theorem there exist a fixed point $u_\infty \in \cA(\cD,\cE, {\bf \mc})$ of $\cT_{\tau}[u]$: 
$$\oT[u_\infty]=u_\infty.$$

Part iv) together with the fact that $\Phi'_k$ and $\Theta'_k$ are positive on $\fR$ implies that $u_\infty'\arrowvert_\fR>0$, and therefore $u_\infty$ is injective on some neighborhood of $I_1$. No conclusion about the injectivity of $u$ can be made on all of $\cD$ since the equality
$$u_\infty= \lambda_\infty^{-1} u_\infty \circ \Tble \circ V_{u_\infty,\tau} \circ  \Tble^{-1}$$
implies that $u'_\infty=0$ whenever $\partial V_{u_\infty,\tau}  \circ  \Tble^{-1} =0$ which might happen at points where $-2 \alpha^{-1}  \Tble^{-1}(z)+u'(z)-\tau'(u(z)) u'(z)=0$ outside of the real line.

We conclude that
$$\eta(z)=u_\infty(T_{b_\infty,\lambda_\infty}(-\sqrt{b_{\infty}-z})), \quad \zeta(z)=u_\infty(T_{b_\infty,\lambda_\infty}(\sqrt{b_{\infty}-z})),$$
are the factorized inverses of a solution $\phi_{\tau}$ of $(\ref{tau_equation})$ on some complex neighborhood of
$$u_\infty(I_1) \supset \left( \max_{(t,s) \in \cS} \Um(l_1;t,s),\min_{(t,s) \in \cS} \um(r_1;t,s) \right) \supset (-1,1)$$.



\section{Properties of $\phi_\tau$}

In this Section we will complete our ``proof'' of the Main Claim.

\noindent
\begin{figure}[t]
  \begin{center}
 \resizebox{70mm}{!}{\includegraphics{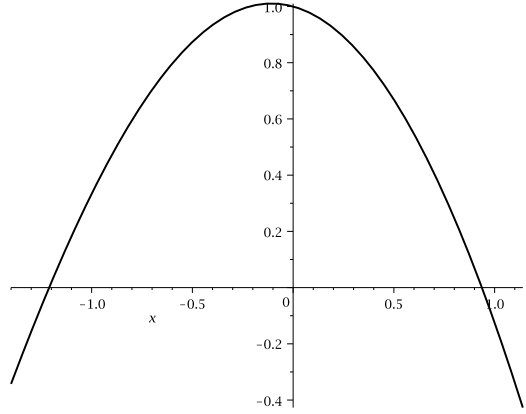}} 
  \caption{ The graph of $\phi_0$}
\label{fig4} 
  \end{center}
\end{figure}

\noindent{\it Demonstration of  the Main Claim $\ref{main_thm}$:}

\noindent For a demonstration of parts $1)$, $2)$, $3)$ and $4)$ see Section $\ref{epsneq0}$ and the Appendix.

\medskip

\noindent {\it Proof of parts  $5)$ and $6)$.} Clearly, since $\phi_\tau(0)=1$ and $\phi_\tau(1)=\lambda/2<0$, there is $0<x_+<1$ such that $\phi_\tau(x_+)=0$. 

Next, consider $\tau=0$. Let $x$ be  in $[x_+ / \lambda,0]$, then $\phi_0(\lambda x)$ assumes all values in $[0,1]$ and $\phi_0(x) \equiv 2 \phi_0(\phi_0(\lambda x))/\lambda -x$ assumes all values in $[(2-x_+)/\lambda,1]$. Since $x_+<1$, we have $(2-x_+)/\lambda <0$ and the equation
$$\phi_0(x)=0$$
has at least one solution $x_- \in [x_+ / \lambda,0]$.

The existence of $x_*$ is obvious. 

We have 
$$\phi_0(1)={2\over \lambda} \phi_0( \phi_0 (\lambda))-1=> {\lambda \over 2} < \phi_0(\phi_0(\lambda)) => \phi_0(\lambda)<1.$$

Furthermore,
$${2\over \lambda } \phi_0(\phi_0(\lambda))={\lambda \over 2}+1 >0 => \phi_0(\phi_0(\lambda))<0 => \phi_0(\lambda)>x_+.$$

Next we demonstrate that $c<0$. 

Since $\phi_0$ is smooth on $I=(-1,1)$, and since $\phi_0(\lambda)<1$, $\phi_0(0)=1$ and $\phi_0(1)<0$,  we must have $c \in (\lambda,1)$. 

Suppose $c>0$. Then $\lambda c \in (\lambda,0)$, and since $\phi_0$ has a single critical point on $I$, $\phi_0$ is strictly decreasing on $I \cap {\fR}_-$ and  $\phi_0(\lambda c)<1$ and $\phi'_0(\lambda c)>0$. 

Now, consider 
$$g(x) \equiv \phi_0(-x)- x={2\over \lambda} \phi_0(\phi_0(-\lambda x)).$$

Clearly, there is an $x_* \in (0,-x_-)$, such that $g(x_*)=0.$

Notice that 
$$g(x_*)=0 => \phi_0(-\lambda x_*)=x_\pm,$$
and since $-\lambda x_* \in (-\lambda x_-,0)$, we have 
\begin{equation}\label{phi_x}
\phi_0(-\lambda x_*)=x_+.
\end{equation}

\noindent
\begin{figure}[t]
  \begin{center}
 \resizebox{70mm}{!}{\includegraphics{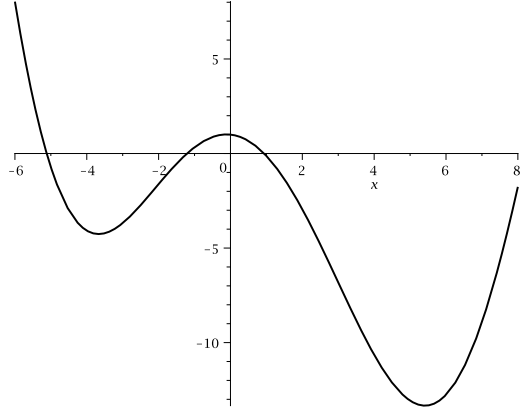}} 
  \caption{ The graph of $\phi_0$ on a larger domain.}
 \label{fig5} 
  \end{center}
\end{figure}

Now suppose that $x_* \le c$, then $-\lambda x_* \in (0,c]$ and 
$$\phi_0(-\lambda x_*)>1,$$ 
which contradicts $(\ref{phi_x})$. Therefore, $x*>c$, and since $g$ is decreasing on $(0,x*)$, we have 
$$g(c)>0 => \phi_0(\lambda c) > \phi_0(-c)+{b-\phi_0(-c) \over 2}>{1 +\lambda \over 2} b +{1-\lambda \over 2}  c >{1 +\lambda \over 2} b c +{1-\lambda \over 2} c >c,$$
where we have used convexity. Finally

$$\phi'_0(c)=2 \phi'_0(\phi_0(\lambda c)) \phi'_0(\lambda c) -1,$$
and since $\lambda c <0 => \phi'_0( \lambda c)>0$, and $\phi_\tau(\lambda c)>c => \phi'_0(\phi_\tau(\lambda c))<0$, we have that 
$$\phi'_0(c)<-1,$$
a contradiction.

{\flushright $\Box$}


\section{Appendix}


\subsection{New a-priori bounds on ${\fR}$}\label{new_bounds}

In this subsection we will use {\it a-priori} bounds on ${\cA}_1({\bf c}))$ to produce quite better bounds on a subset of functions bounded on $(-1,1)$ by a constant.

As before, we denote $(t,s)=(u'(-1/2),u'(0) )$ for a function $u \in  \cA(\cD,\cE, {\bf \mc})$, ${\bf \mc}=(-1/2,0,0,1)$. Recall that $u= \Theta_2 \circ f \circ \Phi_1$, where $f \in \cA_1({\bf c})$, ${\bf c}=(\Phi_1(c_1),\Phi_1(c_2),\Phi_2(c_3),\Phi_2(c_4))$ (note, we will be using the superscript $k$ on functions and numbers, whenever convenient, to avoid double subscripts, these by no means signify raising to a power). Therefore, the following are the derivatives of $f$ at points $c_1$ and $c_2$: 
$$T(t)={t \over \Theta'(c_3) \Phi'_1(-1/2)}, \quad S(s)={s \over \Theta'(c_4) \Phi'_1(0)}.$$

Now, recall that $f'(x)$ is convex, and therefore, using $(\ref{second_der})$,
\begin{eqnarray}\label{2der_bounds}
 & &\min_{x \in [c_1,c_2]} f''(x)  \ge -2 {(c_2-x) T(t)+(x-c_1) S(s) \over (c_2-c_1) (1+c_1)} \equiv m(x,t,s), \\
 & &\max_{x \in [c_1,c_2]} f''(x) \le 2 { (c_2-x) T(t) +(x-c_1) S(s) \over (c_2-c_1) (1-c_2)} \equiv M(x,t,s).
\end{eqnarray} 
Now, fix $t$ and $s$, and consider the function $y(x)=T(t)+\int^x_{c_1} m(z,t,s) d z$. Suppose, that the line $w(x)=S(s)+n(t,s) (x-c_2)$ intersects $(x,y(x))$ at point $x(t,s)$, and $n(t,s)$ is such that the following holds:
\begin{eqnarray} 
c_4-c_3&=&\int^{c_2}_{c_1} \my(z) d z, \\
\my(x)&=&\left\{y(x), \quad c_1 \le x \le x(t,s), \atop w(x), \quad x(t,s) \le x \le c_2. \right.
\end{eqnarray}
First, notice, that any curve $(x,f'(x))$ on $(c_1,c_2)$ with end points $(c_1,t)$ and $(c_2,s)$ can not intersect $(x,y(x))$, and has to intersect $(x,w(x))$ somewhere on $(x(t,s),c_2)$ once ($f'(x)$ is convex), for if it does not then $\int_{c_1}^{c_2} f'(z) d z  \ne 1$. It is also clear that 
\begin{equation}
f(x)  \ge c_3+ \int_{c_1}^x y(z) d z \equiv f_2(x;t,s), \quad  x \in \left[c_1,c_2 \right].
\end{equation}

One can repeat a similar argument for $Y(x)=S(s)+\int^x_{c_2} M(z,t,s) d z$ and $W(x)=T(t)+N(t,s)(x-c_1)$ that intersect at $X(t,s)$ to get 
\begin{eqnarray}
f(x) & \le & c_4-\int_x^{c_2} \mY(z) d z \equiv F_2(x;t,s),  \quad  x \in \left[c_1,c_2 \right], \\ 
\mY(x)&=&\left\{Y(x), \quad  X(t,s) \le x \le c_2, \atop W(x), \quad c_1 \le x \le X(t,s). \right. 
\end{eqnarray}

To obtain an  upper bound on $(-1,c_2)$ and a lower bound on $(c_2,1)$, we recall that the positivity of the Schwarzian derivative for functions in ${\cA}_1({\bf c}))$ together with the positivity of all $f^{(n)}$ for odd $n$ implies that for all $x \in (-1,1)$ 
\begin{equation}\label{S_der}
f'''(x) \ge {3 f''(x)^2 \over 2 f'(x)},
\end{equation}
and consequently,
$$f''(x) \le f''(c_1)+{3 \over 2} \int_{c_1}^x {f''(y)^2 \over  f'(y)} d y,$$
for all $x \in (-1,c_1)$, the equality being realized by the the solution 
$$f'(x)={4 f'(c_1)^3 \over ( -f''(c_1) (x-c_1) +2 f'(c_1) )^2}$$
of equation $(\ref{S_der})$. Therefore, 
$$f(x) \le \int_{c_1}^x {4 T(t)^3 \over ( -f''(c_1) (x-c_1) +2 T(t) )^2}, $$
for all $x \in (-1,c_1)$, the maximum of the right hand side being realized by the maximum admissible $f''(c_1)$ which can be obtained from the condition 
\begin{equation}\label{S_derr}
 {4 T(t)^3 \over ( -f''(c_1) (c_2-c_1) +2 T(t) )^2} = S(s).
\end{equation}
We denote $Z(t,s)$ the solution $f''(c_1)$ of this equation, then 
$$f(x) \le {4 T(t)^3 \over Z(t,s) }\left({1 \over 2 T(t) }- {1 \over 2 T(t) +Z(t,s) (c_1+1) } \right) \equiv F_1(x;t,s),$$
for all $x \in (-1,c_1).$ 

In a similar way
$$  f(x)  \ge {4 S(s)^3 \over X(t,s) }\left({1 \over 2 S(s) +X(t,s) (c_2-1) } - {1 \over 2 S(s) } \right)=f_3(x;t,s), \quad x \in (c_2,1),$$
here $X(t,s)$ solves
\begin{equation}\label{S_derrr}
 {4 S(s)^3 \over ( -X(t,s) (c_1-c_2) +2 S(s) )^2} = T(t).
\end{equation}

Finally, suppose that $ \mm \le f(x) \le \mM$ on the real slice of its domain (this is certainly true if $f  \in  \cA_1({\bf c})$). Consider the line  $(x,S(s)+\mK (x-c_2))$ where $\mK$ is such that 
$$\int_{c_2}^1  S(s)+\mK (x-c_2) d x =\mM-c_4,$$
that is 
$$\mK=2{\mM-c_4 \over (1-c_2)^2}-{S(s) \over 1-c_2}.$$ 
Since $f'(x)$ is convex, the curve $(x,f'(x))$ intersects the line $(x,S(s)+\mK (x-c_2))$ strictly once on $(c_2,1)$.
Convexity of $f'(x)$ implies that 
$$\int_{c_2}^x f'(y) dy<\int_{c_2}^x  S(s)+\mK (y-c_2) d y, \quad x \in (c_2,1),$$
that is
\begin{equation}
f(x) \le c_4+S(s)(x-c_2) +(\mM-c_4-S(s) (1-c_2) ) {(x-c_2)^2 \over (1-c_2)^2 } \equiv F_3(x;t,s), \quad x \in (c_2,1).
\end{equation}

A similar argument on $(-1,c_1)$ demonstrates that 
\begin{equation}
f(x)\ge c_3 -T(t) (c_1-x) + (T(t) (1+c_1)+\mm-c_3){(x-c_1)^2 \over (1+c_1)^2}  \equiv f_1(x;t,s), \quad x \in (-1,c_1).
\end{equation}

Finally, $\mf(x;t,s) \le f(x) \le \mF(x;t,s)$ on $(-1,1)$, where
\begin{equation}\label{f_bounds}
\mf(x;t,s)=\!\left\{\!\begin{array}{cc} f_1(x;t,s),  x \in & \left(-1,c_1 \right) \\  f_2(x;t,s)),  x \in & \left(c_1,c_2 \right) \\  f_3(x;t,s)),  x \in & \left(c_2,1\right) \end{array}  \!\right., \quad    \mF(x;t,s)=\!\left\{\!\begin{array}{cc} F_1(x;t,s),  x \in & \left(-1,c_1\right) \\  F_2(x;t,s)),  x \in & \left(c_1,c_2 \right) \\  F_3(x;t,s)),  x \in & \left(c_2,1\right) \end{array} \! \right.
\end{equation}

Bounds $(\ref{f_bounds})$ transferred to the space $\cA(\cD,\cE;{\bf \mc})$ will be denoted $\um$ and $\Um$:
\begin{eqnarray}
\label{uu_bounds} \um(x;t,s) &\equiv& \Theta_2(\mf_2 (\Phi_1(x);t,s)),\\
 \Um(x;t,s) &\equiv& \Theta_2(\mF_2 (\Phi_1(x);t,s)).
\end{eqnarray}

\subsection{Set of realizable $(u'(-1/2),u'(0))$}\label{ts}

\begin{figure}
 \begin{center}
 \resizebox{50mm}{!}{\includegraphics{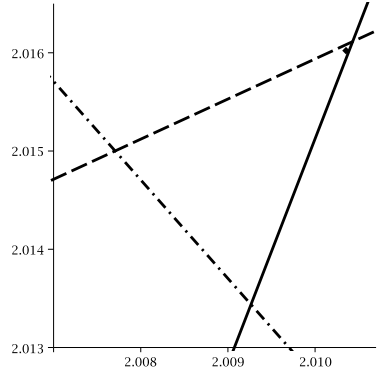}} 
\caption{\it Set $\cS$ bounded by curves $\cZ$ (solid lines), $\cC$ (dash line) and $\cK$ (dash-dot line). The cross mark the pair $(t^*,s^*)$ corresponding to the numerical approximation of the fixed point $u$.}
\label{setS}
\end{center}
\end{figure}

In this subsection we will describe the set $\cS$ of realizable $t=u'(-1/2)$ and $s=u'(0)$ whenever $u \in \cA(\cD, \cE, {\bf  c'})$. We write   $u= \Theta_2 \circ f \circ \Phi_1$,  $f \in \cA_1({\bf c})$ as before.
 
Since $f'(x)<F_1(x;t,s)$ on $(-1,c_1)$ (see Subsection $\ref{new_bounds}$) we have
\begin{equation}\label{s_rs}
 -1+\int_{-1}^{c_1} F_1(x;t,s) d x\le c_3.
\end{equation}
 The relevant (positive) solution $s=s(t)$ of  the equation  $(\ref{s_rs})$ will be denoted  by $\cZ(t)$.

Similarly,
\begin{equation}
 \int_{c_2}^{1} f_3(x;t,s) d x+ c_4 \le 1.
\end{equation}
 The relevant solution $s=s(t)$ of this equation will be denoted by $\cC(t)$.

Finally, the equations
\begin{eqnarray}
\nonumber  \int_{c_1}^{c_2}  {4 T(t)^3 \over ( -Z(t,s) (x-c_1) +2 T(t) )^2}  =c_4-c_3,  \\
\nonumber  \int_{c_2}^{c_1}  {4 S(s)^3 \over ( -X(t,s) (x-c_2) +2 S(t) )^2}  =c_3-c_4, 
\end{eqnarray} 
give two more extremal solution $s=s(t)$ which we will denote $\cG(t)$ and $\cK(t)$.

We have obtained symbolic (and not just numeric) expressions for $\cZ(t)$, $\cC(t)$, $\cG(t)$ and $\cK(t)$ using the Maple software package. The set $\cS$ bounded by these curves is depicted in Fig. $\ref{setS}$.

\subsection{Demonstration of part i) of Claim  $\ref{central_prop}$.}\label{part1}

We will demonstrate existence of a solution $(e,b,\lambda)$ of the equations $(\ref{e_equation})-(\ref{l_equation})$ for the case $\tau \equiv 0$. Specifically, we will show that the set defined by conditions $(\ref{e_bounds})-(\ref{b_bounds})$ is mapped compactly into itself. This, together with the continuity of $\oT$ in $\tau$ at $\tau=0$, implies that there exist $\delta>0$ and $\kappa>0$ such that the same is true for $\tau$'s, holomorphic on $\cE$, whose norm and the norm of whose derivative are bounded  by $\delta$ and $\kappa$ respectively.

To demonstrate $(\ref{e_bounds})$ we introduce a function
$$ \mE(x;\lambda,b)  \equiv -u'(\Tble(x)) {\alpha(b,\lambda) \over 2}.$$
Notice
\begin{eqnarray}
\nonumber \mE(-\gamma(b);\lambda,b)  & \equiv& -u'(\Tble(-\gamma(b))) {\alpha(b,\lambda) \over 2}  = -t {\alpha(b,\lambda) \over 2}, \\
\nonumber \mE(-\beta(b,\lambda); \lambda,b) & \equiv & -u'(\Tble(-\beta(b,\lambda))) {\alpha(b,\lambda) \over 2} =  -s {\alpha(b,\lambda) \over 2}.  
\end{eqnarray}
We verify numerically that $\mE$ maps the interval $(-\beta,-\gamma)$ into itself for all $\cL_-(t,s) \le \lambda \le \cL_+(t,s)$, $\cB_-(\lambda,t,s) \le  b \le \cB_+(\lambda,t,s)$ and $(t,s)\in  \cS$ numerically.

 To show $(\ref{l_bounds})$--$(\ref{b_bounds})$ we consider two functions 
\begin{eqnarray}
\nonumber \mL(\lambda,b) &\equiv& \cT[u](0)=u\left(\Tble\left(\sqrt{b-u\left(\Tble\left(-\sqrt{b-{\lambda \over 2} -{\lambda^2 \over 4}}  \right)\right)} \right) \right), \\
\nonumber  \mB(\lambda,b;e) &\equiv& u \left(\Tble \left(-\sqrt{b-{\lambda \over 2} \left(b-e^2+u(\Tble(e)) \right) } \right) \right),
\end{eqnarray}
and demonstrate that the map
$$(\lambda,b) \mapsto \left(\mL(\lambda,b),\mB(\lambda,b;e) \right)$$
maps  the parallelogram $(\ref{l_bounds})$--$(\ref{b_bounds})$ in the $(\lambda,b)$-plane into itself for all $(t,s) \in \cS$ and all $e$ as in $(\ref{e_bounds})$.

To this end, we first show that $\mL(\cL_+(t,s),b)<\cL_+(t,s)$ for all $\cB_-(\cL_+(t,s),t,s) \le  b \le \cB_+(\cL_+(t,s),t,s)$, and $\mL(\cL_-(t,s),b)>\cL_-(t,s)$ for all  $\cB_-(\cL_-(t,s),t,s) \le  b \le \cB_+(\cL_-(t,s),t,s)$. For this, we define 
\begin{eqnarray}
\nonumber \mL_+(\lambda,b;t,s) \equiv \Um\!\left(\!T_{b,\lambda}\!\left(\!\!\sqrt{b\!-\!\Um\!\left(\!T_{b,\lambda}\!\left(\!-\sqrt{b\!-\!{ \lambda \over 2}\!-\!{ \lambda^2\over 4}}\right)\!;t,s\right)};t,s\!\right) \!\! \right) \\
\nonumber  \mL_-(\lambda,b;t,s) \equiv \um \!\left(\!\!T_{b,\lambda }\!\left(\!\!\sqrt{b\!-\!\um\left(\!T_{b,\lambda}\!\left(\!-\sqrt{b\!-\!{ \lambda \over 2}\!-\!{\lambda^2 \over 4}}\right)\!;t,s\right)};t,s\right) \!\! \right),
\end{eqnarray}
and  verify that, first, 
\begin{equation}\label{l1}
\cL_+(t,s) > \mL_+(\cL_+(t,s),b;t,s)
\end{equation}
for all  $\cB_-(\cL_+(t,s),t,s) \le  b \le \cB_+(\cL_+(t,s),t,s)$ and $(t,s) \in \cS,$
and, second,
\begin{equation}\label{l2}
\cL_-(t,s) < \mL_-(\cL_-(t,s),b;t,s)
\end{equation}
for all  $\cB_-(\cL_-(t,s),t,s) \le  b \le \cB_+(\cL_-(t,s),t,s)$ and $(t,s) \in \cS.$

Inequalities $(\ref{l1})$ and $(\ref{l2})$ have been verified on a computer.

Next, we check that
\begin{eqnarray}
\nonumber \cB_+(\mL(\lambda,\cB_+(\lambda,t,s)),t,s) & > & \mB(\lambda,\cB_+(\lambda,t,s);e), \quad {\rm and} \\
\nonumber \quad \cB_-(\mL(\lambda,\cB_-(\lambda,t,s)),t,s) &< & \mB(\lambda,\cB_-(\lambda,t,s);e),
\end{eqnarray}
for all $\cL_-(t,s) \le \lambda \le \cL_+(t,s)$. We consider
\begin{eqnarray}
\nonumber \mB_+(\lambda,b;e) \equiv \Um\left(T_{b,\lambda}\left(-\sqrt{b-{\lambda \over 2} \left(T^{-1}_{b,\lambda}(e)^2+\Um(e;t,s) \right)}  \right);t,s  \right),\\
\nonumber  \mB_-(\lambda,b;e) \equiv \um\left(T_{b,\lambda}\left(-\sqrt{b-{\lambda \over 2}\left(T^{-1}_{b,\lambda}(e)^2+\um(e;t,s) \right)}  \right);t,s  \right),
\end{eqnarray}
and check the following inequalities numerically:
\begin{eqnarray}
\nonumber \cB_+\left(\mL_-(\lambda,\cB_+(\lambda,t,s)),t,s)  \right) &>& \mB_+(\lambda,\cB_+(\lambda,t,s);e),\\
\nonumber \cB_-(\mL_+(\lambda,\cB_-(\lambda,t,s)),t,s) &<& \mB_-(\lambda,\cB_-(\lambda,t,s);e),
\end{eqnarray}
for all  $\cL_-(t,s) \le \lambda \le \cL_+(t,s)$, $e$ as in $(\ref{e_bounds})$ and $(t,s) \in \cS$.

Finally, since the maps 
$$u \mapsto (\mE(x;\lambda,b),\mL(\lambda,b),\mB(\lambda,d;e)) \quad {\rm and} \quad u \mapsto (\mE(x;\lambda,b),\mL(\lambda,b),\mB(\lambda,d;e))$$
are clearly continuous, we have that there exists at least one fixed point of the map $(\mE,\mL,mB)$ in the set  $(\ref{e_bounds})-(\ref{b_bounds})$ which moves continuously with $u$ and $\tau$.
{\flushright $\Box$}

\subsection{Demonstration of part ii) Claim  $\ref{central_prop}$.}\label{part1}

Differentiate $\cT[u]$ with respect to $x$:
\begin{equation}
\nonumber \cT[u]'(x)= {\alpha^2 \over 2}u' \left( \Tble \left( \Veut( \Tble^{-1}(x)) \right)\right) {u'\left(\Tble \left(-\sqrt{w(x) }\right)\right) \over  4 \Veut(\Tble^{-1}(x) )} {w'(x) \over \sqrt{w(x)}  },
\end{equation}
where  $w$ is the function defined in $(\ref{w_function})$ with $\tau \equiv 0$.

On the real line $ \mw \le  w \le \mW$ and  $\mv \le \Veut \circ \Tble^{-1} \le \mV$ where
\begin{eqnarray}
\nonumber  \mW (x;t,s) &=&{b- {\lambda \over 2} \left(b-\Tble^{-1}(x)^2+  \Um(x;t,s)  \right)},\\
\nonumber  \mw (x;t,s) &=&{b- {\lambda \over 2} \left(b-\Tble^{-1}(x)^2+  \um(x;t,s) \right)},\\
\nonumber \mV (x;ts) &=& \left[{b-\um  \left( \Tble  \left( - \left[ {\mw(x;t,s)} \right]^{1 \over 2}\right);t,s\right)} \right]^{1 \over 2},\\
\nonumber \mv (x;ts) &=& \left[{b - \Um  \left( \Tble  \left( - \left[ {\mW(x;t,s)} \right]^{1 \over 2}\right);t,s\right)} \right]^{1 \over 2}
\end{eqnarray}
are upper and lower bounds on the corresponding functions.

Finally, notice that 

$$
u'(x) \le \Theta_2'(\mF(\Phi_1(x);t,s) ) \mDf(\Phi_1(x);t,s) \Phi_1'(x) \equiv \mDu(x;t,s)
$$
where $\mF$ and $\mf$ are as in $(\ref{f_bounds})$ and 
$$ \mDf(x;t,s) \equiv \eta(x-c_1) \mF(x,t,s) {(1-c_1)  \over  (x-c_1) (1-x) } + \eta(c_1-x)  \mF(x,t,s) {(1+c_1)  \over  (x-c_1) (1+x) } $$
is an upper bound on derivatives on $\cA_1({\bf c})$ that follows from $(\ref{first_der})$ ($\eta$ is the Heaviside function).

Therefore,
\begin{eqnarray}\label{Tu_der}
\nonumber \cT[u]'(x) \le  {\alpha^2 \over 2} \mDu \left( \Tble \left(\mv(x;t,s) \right)\right)\!\!\!\!\!\!\!&\!\!\!\!\!\!\! &\!\!\!\!\!\! {\mDu \left(\Tble \left(-\sqrt{\mW(x;t,s) } \right) \right) \over  4 \mv(x;t,s)} \times \\
&\times &  {2 \alpha^{-1} \Tble^{-1}(x) +(\omega+\sigma x) (1\!+\!\kappa) \over \sqrt{\mw(x;ts)}  }.
\end{eqnarray}

We finally verify on the computer that the right hand side of $(\ref{Tu_der})$ for $\kappa=0$ is strictly less than $\omega+\sigma x$ for all $x \in (0,r_1)$. As in the previous part, this implies existence of $\delta>0$ and $\kappa>0$ such that
$$\oT[u]'(x) < \omega+\sigma x, \quad x \in (0,r_1)$$
whenever $\sup_{z \in \cE}|\tau(z)|<\delta$ and  $\sup_{z \in \cE}|\tau'(z)|<\kappa$.

{\flushright $\Box$}

\subsection{Demonstration of parts iii) and iv) of Claim  $\ref{central_prop}$}\label{part23}

Let $\tau \equiv 0$. We consider a bound on 
\begin{equation}
\cW(\theta) \equiv  w( \partial \cD(\theta)   )
\end{equation}
for $0 \le \theta \le \pi$. Suppose that the boundary of $\cE \cap \fC_+$ is parametrized by $p \in (0,\pi)$. Then, for every fixed $\theta \in [0,\pi]$,  $\cW(\theta)$ is  bounded by the curves 
\begin{eqnarray}
\cW_J(\theta,p)&=&b-{\lambda \over 2} \left(b-T^{-1}_{b,\lambda}(\partial \cD (\theta) )^2 +(p\cE(0) + (1-p) \cE(\pi) ) \right), \quad p \in (0,1),\\
\cW_\cE(\theta,p)&=&b-{\lambda \over 2} \left(b-T^{-1}_{b,\lambda}(\partial \cD (\theta) )^2 +\partial \cE(p) \right), \quad  p \in (0,\pi).
\end{eqnarray}

\begin{figure}[t]
 \begin{center}
\resizebox{70mm}{!}{\includegraphics{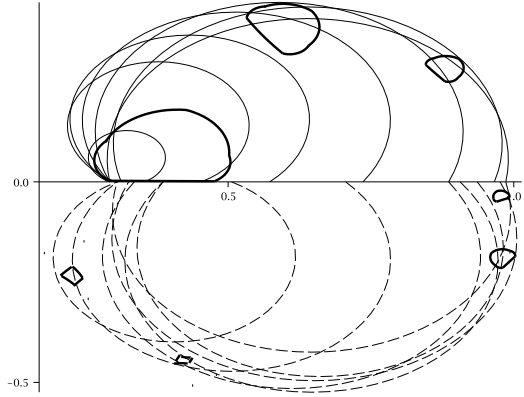}}
\caption{\it Orbit of the set $\cF(\theta)$ for several values of $\theta$. The collection of Poincar\'e neighbourhoods $\cup_k \cD_+(I_k,\theta_k)$ is given in solid lines, the collection $\cup_j \cD_-(I_j,\theta_j)$  --- dashed.}
\label{orbit_fig}
\end{center}
\end{figure}

Denote, $D_5\equiv D(I_5,\vartheta_5)$ (see $(\ref{cond2})$), and let $p \mapsto \partial D_5(p)$ be some parametrization of the boundary of this domain. Then $ \cW(\theta)$ is additionally  bounded by 
\begin{eqnarray}
\nonumber \cW_{D_5} (\theta,p)&=&b-{\lambda \over 2} \left(b-T^{-1}_{b,\lambda}(\partial \cD (\theta) )^2+ \partial D_5(p) \right),\\
\nonumber \cW_{C_\Im}  (\theta,p)&=&b-{\lambda \over 2} \left(b-T^{-1}_{b,\lambda}(\partial \cD (\theta) )^2 + i  C_{\Im} (\theta)\left| \Im(\partial \cD(\theta) ) \right|+p\right),\\
\nonumber \cW_{c_\Im}  (\theta,p)&=&b-{\lambda \over 2} \left(b-T^{-1}_{b,\lambda}(\partial \cD (\theta) )^2 + i  c_{\Im} (\theta)\left| \Im(\partial \cD(\theta) ) \right|+p\right),
\end{eqnarray}
for all $0 \le \theta \le \theta^*$ (see assumption $(\ref{cond1})$ on bounded distortion).

Recall, that  $u= \Theta_2 \circ f \circ \Phi_1$. We first cover the set 
\begin{equation}\label{setF}
\cF(\theta)=\Phi_1\left( T_{b,\lambda} \left(-\sqrt{\cW(\theta)}\right) \right), \quad 0 \le \theta \le \pi,
\end{equation}
 by a collection of  Poincar\'e half-neighbourhoods 
$$\cP=\left(\cup_k \cD_+(I_k,\theta_k) \right) \cup \left(\cup_j \cD_-(I_j,\theta_j)\right),$$
for some appropriately chosen $I_k=(l_k,r_k)$, $I_j=(l_j,r_j)$ and $\theta_k$, $\theta_j$ (see Fig. \ref{orbit_fig}), then according to Lemma $\ref{Epstein_lemma}$, the set $f(\cF(\theta))$, $0 \le \theta \le \pi$ is contained in 

$$\tilde{\cP}(t,s)=\left( \cup_k \cD_+(\tilde{I}_k,\theta_k) \right) \cup \left(\cup_j \cD_-(\tilde{I}_j,\theta_j) \right),$$
where 
$$\tilde{I}_m=(\mf(l_m;t,s),\mF(r_m;t,s))$$

Set $\cV(t,s) \equiv \Theta_2(\tilde{\cP}(t,s))$.  We construct the set 
$$\cM(t,s)=-{\rm sign}\left(\Im\left(b-\cV(t,s) \right)  \right)\sqrt{ b-\cV(t,s)},$$
which is a bound on $V_{u,0}(\Tble^{-1}(\cD_+))$, and verify that it is contained compactly in $\Tble^{-1}(\cD_-)$.

\begin{figure}[t]
 \begin{center}
\resizebox{75mm}{!}{\includegraphics{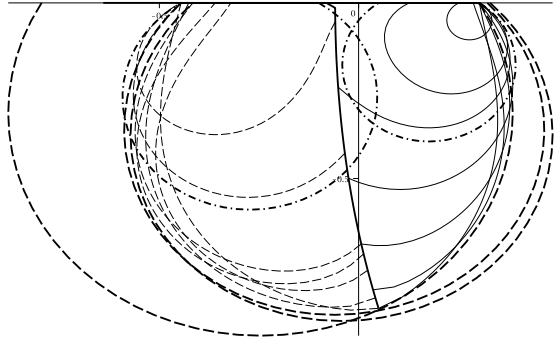}}
\caption{\it Set $\cN$ (thin  lines) covered by a collection $\cH=\cH_1 \cup \cH_2$ of Poincar\'e neighbourhoods:  $\cH_1$ (thick dash lines) -- intersection of three neighbourhoods, $\cH_2$ (thick dash-dot lines) -- union of two. The thick line in the middle is the  image of $\fR \pm i 0$ under the map $z \mapsto \Phi_1\left( \Tble \left(  -{\rm sign}\left(\Im\left[b-u(z) \right] \right)\sqrt{ b-u(z)} \right) \right)$.}
\label{cover_1}
\end{center}
\end{figure}

Notice, that 
$$\bar{H} \equiv (b-\cV(t,s) ) \cup \fR_-$$
is non-empty. Therefore $\sqrt{}$ is not defined on $\bar{H}$, and neither is $V_{u,0}$ on$$H \equiv q^{-1}(\bar{H})), \quad q\equiv b-u\left(\Tble\left(-\sqrt{w(\Tble(z))}  \right)\right).$$
It is easily checked however that it is continuous  across $H$, and holomorphic in both components of $\Tble^{-1}(\cD_+) \setminus H$. Therefore, by Morera's theorem, it is holomorphic in all of $\Tble^{-1}(\cD_+)$. Analyticity on  $\Tble^{-1}(\cD_-)$ follows in a similar way.

We next construct the set 
$$\cN(t,s)=\Phi_1 \left( \Tble \left(\cM(t,s);t,s \right)\right)$$
and cover it with another collection of Poincar\'e  half-neighbourhoods (see  Fig. \ref{cover_1}):
\begin{eqnarray}
\nonumber \cH_1&=& \cap_n \cD_+(K_n,\alpha_n), \quad K_n=(k_n,j_n), \\
\nonumber \cH_2&=& \cup_n \cD_+(J_n,\phi_n), \quad J_n=(m_n,p_n),\\
\nonumber \cH&=&\cH_1 \cup \cH_2.
\end{eqnarray} 

Set
\begin{eqnarray}
\nonumber \tilde{\cH}_1(t,s)&=& \cap_n \cD_+(\tilde{K}_n,\alpha_n), \quad \tilde{K}_n=(\mf(k_n;t,s),\mF(j_n;t,s)),\\
\nonumber \tilde{\cH}_2(t,s)&=& \cup_n \cD_+(\tilde{J}_n,\phi_n),   \quad \tilde{J}_n=(\mf(m_n;t,s),\mF(p_n;t,s)),\\
\nonumber \tilde{\cH}(t,s)&=&\tilde{\cH}_1(t,s) \cup \tilde{\cH}_2(t,s).
\end{eqnarray} 

\begin{figure}[t]
 \begin{center}
\resizebox{60mm}{!}{\includegraphics{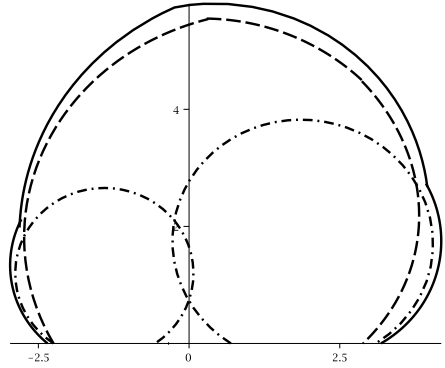}}
\caption{\it Set $\cX$, the union of  $\lambda^{-1} \Theta_2  ( \tilde{\cH}_1)$ (dash lines) and   $\lambda^{-1} \Theta_2  ( \tilde{\cH}_2)$ (dash-dot lines),  is contained in $\cE$ (solid lines).}
\label{cover_2}
\end{center}
\end{figure}

Finally, the set
$$\cX(t,s)=\lambda^{-1} \Theta_2  ( \tilde{\cH}(t,s))$$
is verified numerically to be contained  compactly in $\cE$ (see Fig. $\ref{cover_2}$) for all $(t,s) \in \cS$. This shows that $\cT[u]$ is in $\cA(\cD,\cE,{\bf \mc})$ whenever $u \in \cA(\cD,\cE,{\bf \mc})$.

The claims of parts iii) and iv) in case of non-zero $\tau$ follow from the fact that all set containments verified in the case of $\tau \equiv 0$ were compact.

{\flushright $\Box$}

\subsection{Demonstration of part v) of Claim  $\ref{central_prop}$:    Invariance of bounded  boundary distortion} \label{part5}

As can be seen  from the previous proofs, to demonstrate the invariance of the set $\cA(\cD,\cE,{\bf \mc})$, one needs to  assume a bound on the distortion of $\partial {\cD}$ by $u$. We, therefore, need to reproduce this  bound for ${\cT}[u]$. 

Suppose $\theta \mapsto \partial \cD(\theta)$ is a parametrization of $\partial \cD$  as described in subsection $\ref{part23}$.  Recall, that we have made the following assumption (see $(\ref{cond1})$): 
$$u(\partial \cD(\theta) ) \in \cO(\theta) \equiv \left( \cE  \setminus D(I_5,\theta_5) \right)  \cap \{z \in \field{C}: c_{\Im}(\theta) |\Im(\partial \cD(\theta) )|  < \Im(z) < C_{\Im}(\theta) |\Im(\partial \cD(\theta) )|  \}.$$
 for $0 \le \theta \le \theta^*$.

\begin{figure}
 \begin{center}
\begin{tabular}{c c c}
 \resizebox{49mm}{!}{\includegraphics{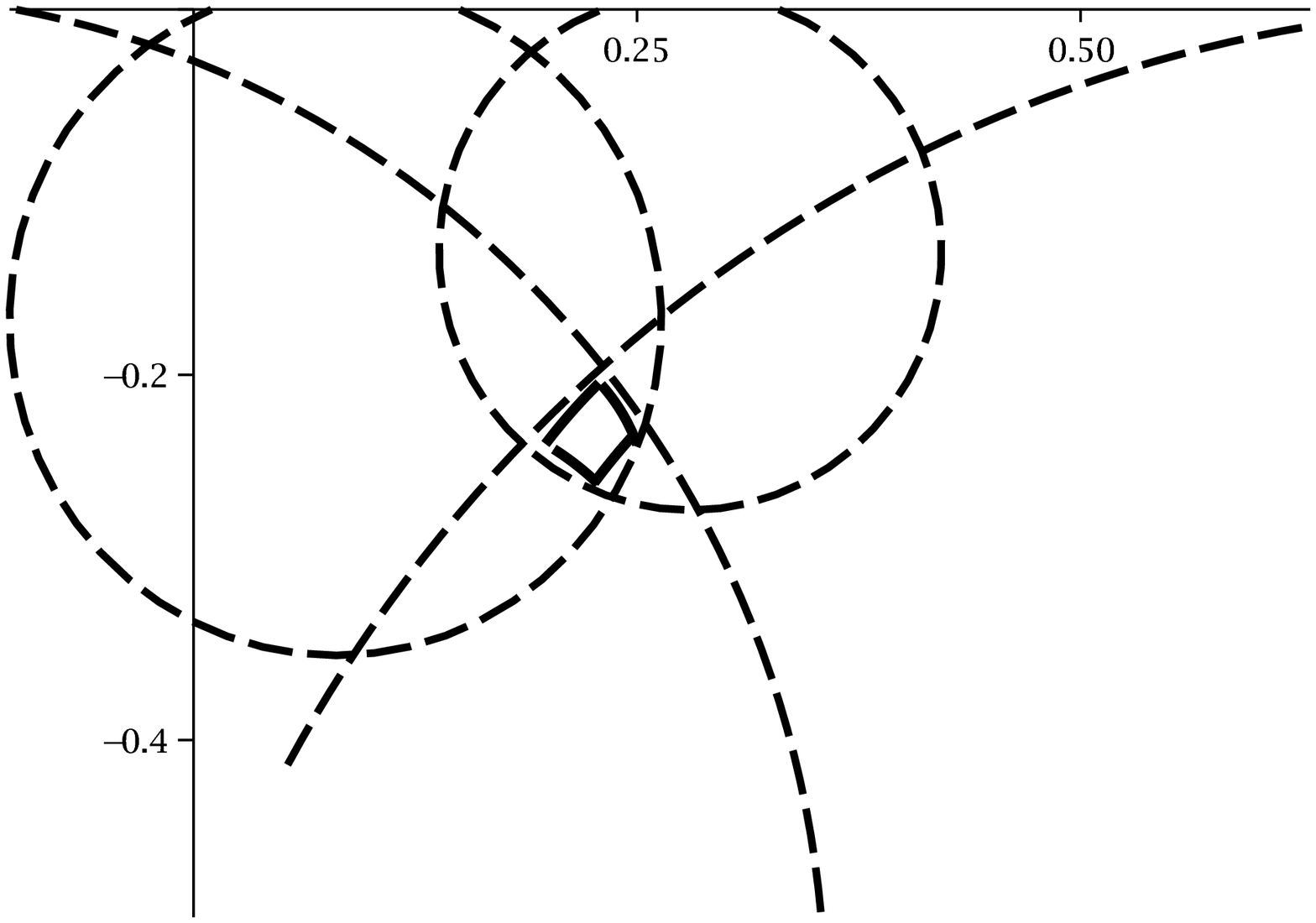}} & \resizebox{47mm}{!}{\includegraphics{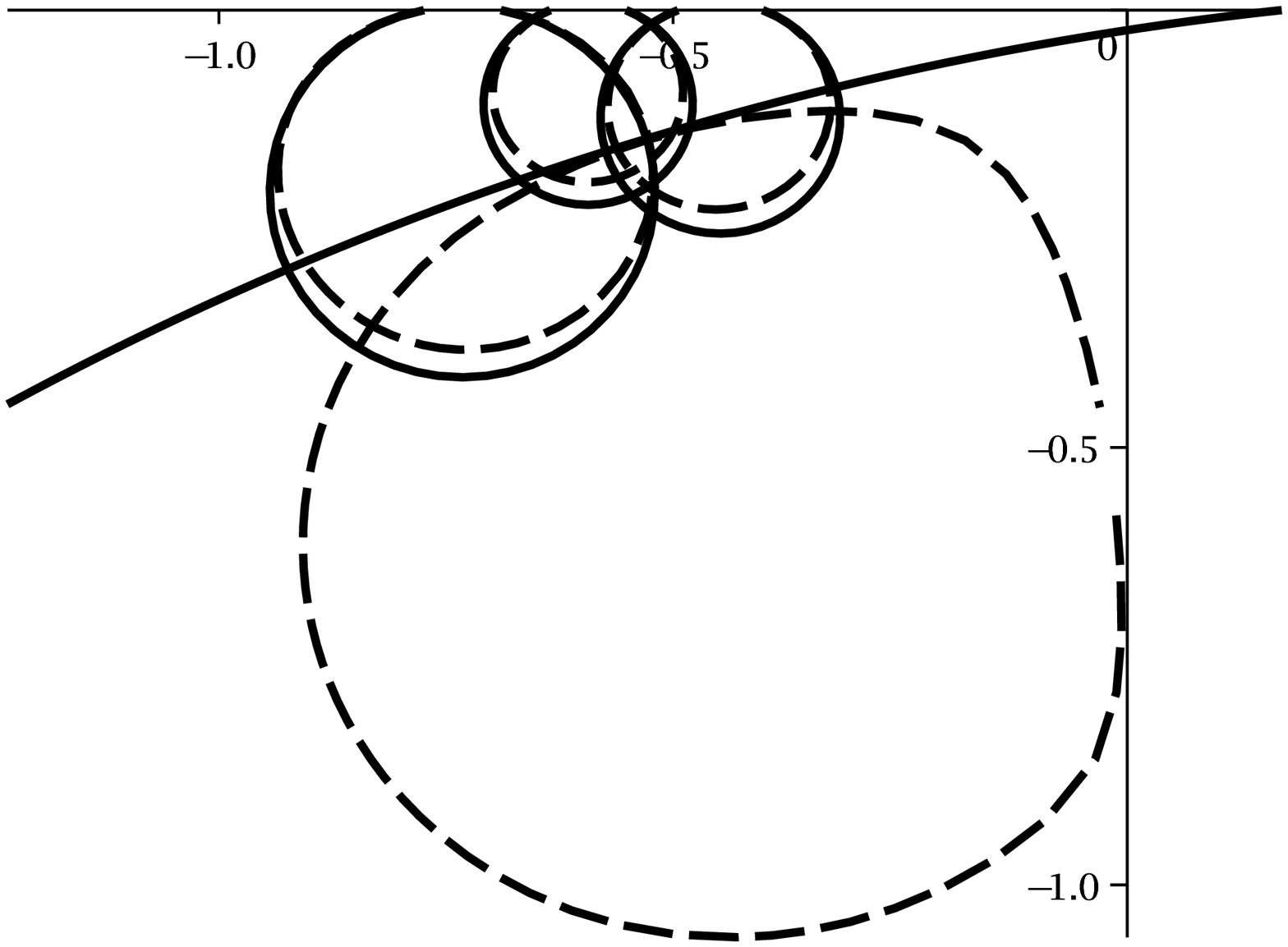}} & \resizebox{49mm}{!}{\includegraphics{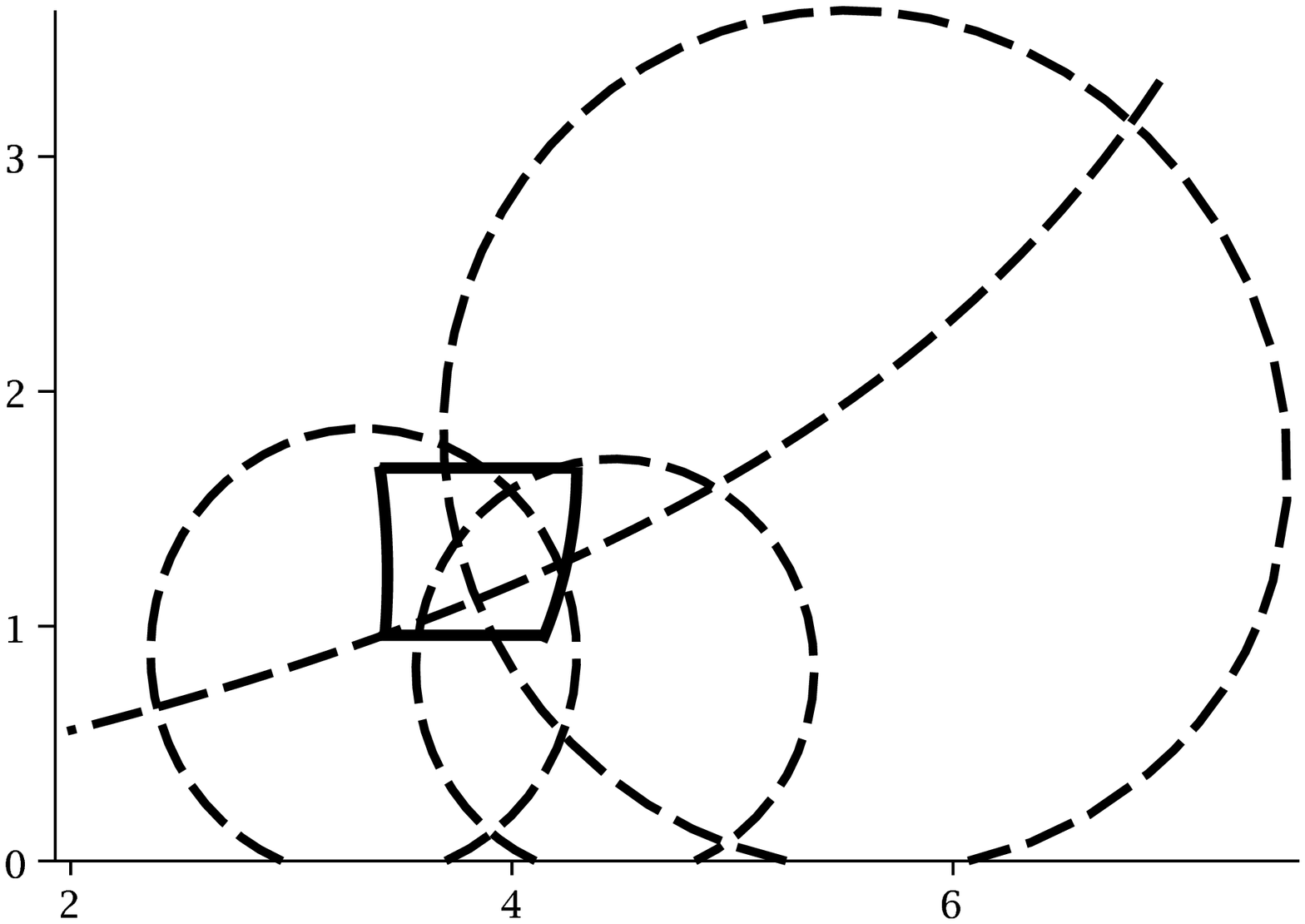}} \\
a) & b) & c
\end{tabular} 
\caption{\it a). Set $\cF(\theta)$ (solid lines), $\theta=0.3$, is contained in the intersection of four Poincar\'e neighbourhoods (dash lines); b). Set $\cQ(\theta)$, intersection of the sets given by dash lines, is contained in the intersection of four Poincar\'e neighbourhoods (solid lines); c)  Set $\cG(\theta)$, the intersection of four sets given by the dash lines, is inside the set $\cO(\theta)$ (solid lines).}
\label{bounds}
\end{center}
\end{figure}

To show that these assumptions are reproduced for $\cT[u]$, we first construct four continuous families of Poincar\'e half-neighbourhoods 
$$D_-(L_i(\theta;b,\lambda),\vartheta_i(\theta;b,\lambda)),  \quad L_i(\theta;b,\lambda)=(w_i(\theta;b,\lambda),v_i(\theta;b,\lambda)) \Subset (-1,1), \quad i=1..4,$$
(below, we will suppress dependence on $b$ and $\lambda$) that contain the set 
$$\cF(\theta) \equiv \Phi_1\left(\Tble \left(-\sqrt{w(\cO(\theta))  } \right)   \right)$$
(see $(\ref{setF})$) in their intersection (Fig $\ref{bounds}$ a)).  Set
$$\cK(\theta)=\Theta_2 \left( \cap_{i=1..4}  D_-(\tilde{L_i}(\theta),\vartheta_i(\theta)) \right), \quad \tilde{L_i}(\theta)=(\mf(w_i(\theta);t,s),\mF(v_i(\theta);t,s)).$$

By Schwarz Lemma,
$$u \left(  T_{b,\lambda} \left(-\sqrt{\cW(\theta)}\right) \right) \subset \cK(\theta).$$
Next, we set 
$$\cQ(\theta)=\Phi_1 \left( T_{b,\lambda} \left( -{\rm sign} \left(\Im \left(b-  \cK(\theta) \right)  \right) \sqrt{ b- \cK(\theta) } \right) \right),$$
and construct four more families of Poincar\'e neighbourhoods  
$$D-(M_i(\theta),\varphi_i(\theta)), \quad i=1..4, \quad M_i=(h_i(\theta),q_i(\theta)) \Subset (-1,1),$$
that contain $\cQ(\theta)$ in their intersection (see Fig $\ref{bounds}$ b)). Again, by Schwarz Lemma,
$$\cT[u](\partial \cD(\theta)) \subset \lambda^{-1} \Theta_2 \left( \cap_{i=1..4} D_-((\mf(h_i(\theta);t,s),\mF(q_i(\theta);t,s)),\varphi_i(\theta)) \right) \equiv \cG(\theta)$$
(we suppress dependence on $t$ and $s$ in $\cG$).

Finally, we verify numerically that 
$$\cG(\theta) \subset \cO(\theta)$$
 for all $0 \le \theta \le \theta^*$, $\cL_-(t,s) \le \lambda \le \cL_+(t,s)$, $\cB_-(\lambda, t,s) \le b \le \cB_+(\lambda, t,s)$ and $(t,s) \in \cS$.

$\Box$


\begin{thebibliography}{****} 

\bibitem[AK]{AK} J. J. Abad, H. Koch, Renormalization and periodic orbits for Hamiltonian flows, {\it Comm.  Math.  Phys. } { \bf 212} (2000) \# 2 371--394.

\bibitem[AKW]{AKW} J. J. Abad, H. Koch and P. Wittwer,  A renormalization group for Hamiltonians: numerical results, {\it Nonlinearity} { \bf 11} (1998) 1185--1194.

\bibitem[BCGG]{BCGG} G. Benettin et al, Universal properties in conservative dynamical systems, {\it Lettere al Nuovo Cimento} {\bf 28} (1980) 1--4.


\bibitem[Bou]{Bou} T. Bountis,  Period doubling bifurcations and universality in conservative Systems, {\it Physica} { \bf 3D} (1981) 577--589.

\bibitem[dCLM]{dCLM} A. de Carvalho, M. Lyubich, M. Martens,  Renormalization in the H\'enon family, I: Universality but non-rigidity, {\it Nonlinearity} { \bf 11} (1998) 1185--1194.

\bibitem[CEK1]{CEK1} P. Collet, J.-P. Eckmann and H. Koch, Period doubling bifurcations for families of maps on ${\fR}^n$, {\it J. Stat. Phys.} { \bf 3D} (1980).

\bibitem[CEK2]{CEK2} P. Collet, J.-P. Eckmann and H. Koch, On universality for area-preserving maps of the plane , {\it Physica} { \bf 3D} (1981) 457--467.

\bibitem[DP]{DP} B. Derrida, Y. Pomeau, Feigenbaum's ratios of two dimensional  area preserving maps, {\it Phys. Lett.} {\bf A80} (1980) 217--219.


\bibitem[EKW1]{EKW1} J.-P. Eckmann, H. Koch and P. Wittwer, Existence of a fixed point of the doubling transformation for area-preserving maps of the plane, {\it Phys. Rev. A} { \bf 26} (1982) \# 1  720--722.

\bibitem[EKW2]{EKW2} J.-P. Eckmann, H. Koch and P. Wittwer, A Computer-Assisted Proof of Universality for Area-Preserving Maps, {\it Memoirs of the American Mathematical Society} {\bf  47} (1984), 1--121.


\bibitem[Eps1]{Eps1} H. Epstein, New proofs of the existence of the Feigenbaum  functions, {\it  Commun. Math. Phys.} {\bf 106} (1986) 395--426.


\bibitem[Eps2]{Eps2} H. Epstein, Fixed points of composition operators II, {\it  Nonlinearity } {\bf 2} (1989) 305--310.


\bibitem[ED]{ED} D. F. Escande, F. Doveil, Renormalization method for computing the threshold of the large scale stochastic instability in two degree of freedom Hamiltonian systems, {\it  J .Stat. Phys. } {\bf 26} (1981) 257--284.  

\bibitem[dF1]{dF1} E. de Faria, {\it Proof of universality for critical circle mappings}, Thesis, CUNY, 1992.

\bibitem[dF2]{dF2} E. de Faria, Asymptotic rigidity of scaling ratios for critical circle mappings, {\it Ergodic Theory Dynam. Systems} {\bf 19}(1999), no. 4, 995--1035.

\bibitem[Fei1]{Fei1} M. J. Feigenbaum,  Quantitative universality for a class of nonlinear transformations, {\it J. Stat. Phys.} {\bf 19} (1978) 25--52.
\bibitem[Fei2]{Fei2} M. J. Feigenbaum, Universal metric properties of non-linear transformations, {\it J. Stat. Phys.} {\bf 21} (1979) 669--706.

\bibitem[Gai1]{Gai1} D. Gaidashev, Renormalization of isoenergetically degenerate Hamiltonian flows and associated bifurcations of invariant tori, {\it Discrete Contin. Dyn. Syst.} {\bf 13}(2005), no. 1, 63--102.

\bibitem[Gai3]{Gai3}  D. Gaidashev,  Computer-assisted bounds on the solution of a Beltrami equation and applications to renormalization, e-print math.DS/0510472 at Arxiv.org.

\bibitem[GK]{GK} D. Gaidashev, H. Koch, Renormalization and shearless invariant tori: numerical results, {\it Nonlinearity} {\bf 17}(2004), no. 5, 1713--1722.

\bibitem[GaiYa]{GaiYa}D. Gaidashev, M. Yampolsky, Cylinder renormalization of Siegel disks, {\it Exp. Math.} {\bf 16:2} (2007).

\bibitem[Hel]{Hel} R. H. G. Helleman, Self-generated chaotic behavior in nonlinear mechanics, in "Fundamental problems in statistical mechanics", Ed. by E. G. D. Cohen, North-Holland, Amsterdam, p.165, (1980).


\bibitem[KLDM]{KLDM}  K. Khanin, J. Lopes Dias, J. Marklof,  Multidimensional continued fractions, dynamic renormalization and KAM theory,  {\it Comm. Math. Phys.},  {\bf 270}  (2007),  no. 1, 197--231.


\bibitem[Koch1]{Koch1}H. Koch, On the renormalization of Hamiltonian flows, and critical invariant tori, {\it Discrete Contin.  Dyn.  Syst.}  {\bf 8} (2002),  633--646.
\bibitem[Koch2]{Koch2}H. Koch, A renormalization group fixed point associated with the breakup of golden invariant tori, {\it  Discrete Contin. Dyn. Syst.}  {\bf 11}  (2004),  no. 4, 881--909.
\bibitem[Koch3]{Koch3} H. Koch, {\it Existence of critical invariant tori}, preprint {\tt mp\_arc 04--210} (2004)
\bibitem[Koci\'c]{Kocic} S. Koci\'c, Renormalization of Hamiltonians for Diophantine frequency vectors and KAM tori,{\it Nonlinearity} {\bf 18} (2005) 2513--2544.

\bibitem[Lyu]{Lyu} M. Lyubich, Feigenbaum-Coullet-Tresser universality and Milnor's  hairness conjecture, {\it Annals of Mathematics} {\bf 149} (1999) 319--420.

\bibitem[LY]{LY} M. Lyubich, M. Yampolsky, Dynamics of quadratic polynomials: complex  bounds for real maps, {\it Ann. Ins. Fourier} {\bf 47} 4 (1997) 1219-1255.


\bibitem[McK1]{McK1} R. S. MacKay, Renormalization approach to invariant circles in area-preserving maps, {\it Physica} {\bf D7} (1983) 283--300. 
\bibitem[McK2]{McK2} R. S. MacKay, Renormalisation in area preserving maps, Thesis, Prin\-ce\-ton (1982). World Scientific, London (1993).
\bibitem[ME]{ME} A. Mehr and D.F. Escande, {\it  Physica }{\bf D13}  (1984) 302.


\bibitem[McM]{McM} C. McMullen, Self-similarity of Siegel disks and Hausdorff dimension of  Julia sets, {\it Acta Math.} {\bf 180}(1998), 247-292.

\bibitem[Shen]{Shen} S. J. Shenker, L. P. Kadanoff, Critical behaviour of KAM surfaces. I Empirical results, {\it  J. Stat.  Phys. } {\bf 27} (1982) 631--656. 

\bibitem[Spa]{Spa} C. Sparrow, The Lorenz equations : bifurcations, chaos, and strange attractors, New York ; Berlin : Springer (1982).

\bibitem[Sul]{Sul} D. Sullivan, Bounds, quadratic differentials and renormalization conjectures, in: Mathematics into the Twenty-first Century, AMS Centennial Publications, Vol. II, Amer. Math. Soc., Providence, R.I. (1992) 417-466.

\bibitem[TC]{TC} C. Tresser and P. Coullet, It\'erations d'endomorphismes et groupe de renormalisation, {\it  C. R. Acad.  Sci. Paris} {\bf 287A}(1978), 577--580.


\bibitem[Ya1]{Ya1} M. Yampolsky, Hyperbolicity of renormalization of critical circle maps,
{\it  Publ. Math. Inst. Hautes Etudes Sci.} {\bf 96}(2002), 1--41.

\bibitem[Ya2]{Ya2} M. Yampolsky, Renormalization horseshoe for critical circle maps,
{\it  Commun. Math. Physics} {\bf 240}(2003), 75--96.

\bibitem[Ya3]{Ya3} M. Yampolsky,  Siegel disks and renormalization fixed points, e-print math.DS/0602678 at Arxiv.org

\end{thebibliography}
\end{document}